\documentclass[11pt]{amsart}

\usepackage{amssymb}
\usepackage{amsmath}
\usepackage{enumerate}
\usepackage{amsfonts}
\usepackage{color}
\usepackage{latexsym}
\usepackage{mathtools}
\usepackage{amstext,amssymb,amsopn,amsthm,mathrsfs,wasysym,dsfont,comment,bigints,tensor,upgreek} 
\usepackage{longtable,tabularx,booktabs}

\usepackage{float}

\usepackage{textcomp}

\usepackage[dvipsnames]{xcolor}
\usepackage{tikz}
\usetikzlibrary{patterns,arrows}
\usepackage{cite}
\usepackage[
colorlinks=true,
linkcolor=magenta,
citecolor=magenta,
urlcolor=magenta,
pagebackref,
]{hyperref}

\numberwithin{equation}{section}

\usepackage{libertine}
\usepackage[libertine]{newtxmath}

\def\L{\mathcal{L}}
\def\grad{\nabla_{\!\H}}
\usepackage[bb=boondox]{mathalfa}
\def\I{\mathbb{I}}
\def\0{\mathbb{0}}
\newcommand{\jap}[1]{\left\langle#1\right\rangle}
\def\J{\mathcal{J}}
\def\kor{\rho}
\def\Kor{\left|(z,t)\right|_{\H}}
\def\u{\mathfrak{u}}
\def\f{\mathfrak{f}}
\def\wgrad{\widetilde{\nabla}_{\!\H}}
\def\gl{\upomega}
\def\e{\varepsilon}
\def\k{\kappa}
\def\l{\lambda}

\def\C{\mathbb{C}}
\def\R{\mathbb{R}}

\def\H{\mathbb{H}}

\def\o{\overline}

\def\til{~}
\definecolor{myyellow}{rgb}{0.9290 0.6940 0.1250}
\definecolor{myorange}{rgb}{0.8500 0.3250 0.0980}
\definecolor{myred}{rgb}{0.6350 0.0780 0.1840}
\definecolor{mygreen}{rgb}{0.4660 0.6740 0.1880}
\definecolor{mycyan}{rgb}{0.3010 0.7450 0.9330}
\definecolor{myblue}{rgb}{0 0.4470 0.7410}
\definecolor{myred}{rgb}{0.4940 0.1840 0.5560}
\newcommand{\norm}[1]{\left\lVert#1\right\rVert}

\DeclareMathOperator*{\sgn}{sgn}

\newtheorem{theorem}{Theorem}

\newtheorem{lemma}{Lemma}

\theoremstyle{definition}

\newtheorem{remark}{Remark}[section]

\newcommand{\im}{\operatorname{Im}}
\newcommand{\re}{\operatorname{Re}}

\allowdisplaybreaks

\author[L.\ Fanelli]{Luca Fanelli}
\address[Luca Fanelli]{
Ikerbasque, Basque Foundation for Science,
48011 Bilbao, Spain,
\newline \phantom{\quad} \&
Universidad del Pa\'is Vasco / Euskal Herriko Unibertsitatea,
48080 Bilbao, Spain}
\email{\href{mailto:luca.fanelli@ehu.es}{luca.fanelli@ehu.es}}

\author[H.\ Mizutani]{Haruya Mizutani}
\address[Haruya Mizutani]{
Department of Mathematics, Graduate School of Science, Osaka University, 
\newline \phantom{\quad} Toyonaka, 560-0043 Osaka, Japan}
\email{\href{mailto:haruya@math.sci.osaka-u.ac.jp}{haruya@math.sci.osaka-u.ac.jp}}

\author[L.\ Roncal]{Luz Roncal}
\address[Luz Roncal]{
BCAM -- Basque Center for Applied Mathematics,
48009 Bilbao, Spain,
\newline \phantom{\quad} \&
Ikerbasque, Basque Foundation for Science,
48011 Bilbao, Spain,
\newline \phantom{\quad} \&
Universidad del Pa\'is Vasco / Euskal Herriko Unibertsitatea,
48080 Bilbao, Spain}
\email{\href{mailto:lroncal@bcamath.org}{lroncal@bcamath.org}}

\author[N. M. Schiavone]{Nico Michele Schiavone}
\address[Nico M. Schiavone]{
	BCAM -- Basque Center for Applied Mathematics,
	48009 Bilbao, Spain}
\email{\href{mailto:nschiavone@bcamath.org}{nschiavone@bcamath.org}}


\makeatletter
\@namedef{subjclassname@2020}{%
  \textup{2020} Mathematics Subject Classification}
\makeatother


\date{\today}

\begin{document}

\title[Non-existence of radial eigenfunctions for the perturbed Heisenberg sublaplacian]{Non-existence of radial eigenfunctions \\ for the perturbed Heisenberg sublaplacian}

\begin{abstract}
	We prove uniform resolvent estimates in weighted $L^2$-spaces for radial solutions of the sublaplacian $\L$ on the Heisenberg group $\H^d$. The proofs are based on the multipliers methods, and strongly rely on the use of suitable multipliers and of the associated Hardy inequalities. The constants in our inequalities are explicit and depend only on the dimension $d$. As application of the method, we obtain some suitable smallness and repulsivity conditions on a complex radial potential $V$ on $\H^d$ such that $\L+V$ has no radial eigenfunctions.
\end{abstract}

\maketitle





\section{Introduction}

In this manuscript, we are interested in the following inhomogeneous Helmholtz equation for the sublaplacian $\L$ in the Heisenberg group $\H^d$
\begin{equation}\label{eq:prob0}
-\L u + \lambda u = f
\end{equation}
where $d\ge1$, $u,f \colon \H^d \to \C$, and $\lambda\in\C$. Merely assuming
$f \in H^{-1}(\H^d)$, we can make sense of the weak formulation in $H^1(\H^d)$ of equation \eqref{eq:prob0}:
we say that $u$ is a \emph{solution} of \eqref{eq:prob0} if $u\in H^1(\H^d)$ and
\begin{equation}\label{eq:solution}
	- \langle \nabla_{\!\H} v,  \nabla_{\!\H} u \rangle
	+
	\lambda
	\langle v, u \rangle
	=
	\langle v, f \rangle
\end{equation}
for all $v \in H^1(\H^d)$, see Section \ref{sec:main_results} for the notation. 
With a customary abuse of notation, in the left-hand side $\langle \cdot, \cdot \rangle$ denotes the inner product on $L^2(\H^d)$, and in the right-hand side the duality pairing between $H^1(\H^d)$ and its dual\til$H^{-1}(\H^d)$.

In the Euclidean setting, equation \eqref{eq:prob0} reads as
\begin{equation}\label{eq:euclidean}
\Delta u+\lambda u=f
\end{equation}
with $u,f:\R^d\to \C$, and $\lambda\in\C$.
The above equation is probably the strongest link between Spectral Theory, Fourier Analysis and Partial Differential Equations. When $f=Vu$, with \mbox{$V:\R^d\to\C$}, this is the eigenvalue equation associated to the operator $-\Delta+V$, which is possibly non self-adjoint, since $V$ is complex-valued. Recently, due to the introduction of quantum mechanical models described by non-self-adjoint Hamiltonians (see e.g. \cite{BGSZ, BB, SGY}), an increasing interest of the mathematical community has been devoted to the associated spectral analysis (we mention, among others, \cite{CFK, CFK23, DFKS, DFS, FKV, FKV2, Frank, FrankIII, FS17, MS, S} and the references therein). Some basic questions about the spectrum of these Hamiltonians are quite difficult to answer, since many classical tools such as the Spectral Theorem and the variational methods are not available in the non-self-adjoint setting. Tipically, the most challenging problems concern the point spectrum: for instance, when an external potential is present, it is desirable to determine a suitable notion of \textit{smallness} such that the potential is not able to create eigenvalues.
Recently, in the non-self-adjoint context, two particular strategies have proven to be quite useful at this aim: on the one hand, the combination of suitable resolvent estimates with the Birman--Schwinger principle; on the other hand, the so-called method of multipliers.

\subsection{Resolvent estimates and Birman--Schwinger principle}

After the pioneer results by Simon in \cite{S}, related with symmetric Hamiltonians, a breakthrough has been provided in the seminal work by Frank \cite{Frank}, who introduced the so-called $L^p-L^q$ uniform resolvent estimates, a classical tool from Fourier analysis and PDEs which mainly refers to the celebrated paper \cite{KRS} by Kenig, Ruiz and Sogge, and its further generalizations \cite{Gut,RXZ,KL}. 
%
%
Combing these estimates with the Birman--Schwinger principle (see \cite{HK} for a comprehensive picture of the principle), Frank showed that the point spectrum of $-\Delta+V$, with $V$ possibly complex-valued, is empty provided that $\|V\|_{L^{d/2}}$ is sufficiently small. 
The proof of the uniform resolvent estimates strongly relies on the Fourier representation of solutions to \eqref{eq:euclidean}, which in the setting of the Heisenberg group is somehow tricky. Indeed, the group Fourier transform of a function  $f\in L^1(\H^d)$ is the operator-valued function defined, for each $\lambda\in \R^*$, by
\begin{equation*}
	\widehat{f}(\lambda) := \pi_{\lambda}(f)= \int_{\H^d}f(z,t)\pi_{\lambda}(z,w)\,dz\,dt,
\end{equation*}
where $\pi_{\lambda}$, for each $\lambda\in\R^*=\R\setminus\{0\}$, is an irreducible  unitary representation of $\H^d$ realized on $ L^2(\R^d)$ with action given by
$$
\pi_\lambda(z,t)\varphi(\xi) = e^{i\lambda t} e^{i(x\cdot \xi+\frac12 x\cdot y)}\varphi(\xi+y) 
$$
for $\varphi \in L^2(\R^d)$ and $ z = x+iy$. 
As a fact, up to our knowledge, uniform estimates are not available in the setting of the Heisenberg group for equation \eqref{eq:prob0}. We refer the reader to the renowned papers by  Folland and Stein \cite{FS}, Jerison and Lee \cite{JL}, and Frank and Lieb \cite{FL} for the related Hardy--Littlewood--Sobolev estimates at the zero-energy $\lambda=0$.

\subsection{Kato--Yajima estimates and multipliers method.}

In the attempt to reach in the Heisenberg context results on resolvent estimates and spectral stability akin to their Euclidean counterparts, we will look for uniform resolvent estimates which can be proved by real analytical techniques. This is the case of weighted $L^2$-estimates, for which a huge amount of results is available. Among them, in \cite{KY} Kato and Yajima proved, along with other inequalities, that
\begin{equation}\label{eq:KY}
\big\||x|^{-1}u\big\|_{L^2(\R^d)} \leq C\||x|f\|_{L^2(\R^d)},
\qquad
\lambda\in\C\setminus[0,+\infty), \quad d\ge 3,
\end{equation}
for some $C>0$ independent on $\lambda$, where $u$ solves \eqref{eq:euclidean}. Notice that the zero-energy version of~\eqref{eq:KY} is
\begin{equation}\label{eq:KYequiv}
	\big\||x|^{-1}(-\Delta)^{-1}f\big\|_{L^2(\R^d)}\leq C\||x|f\|_{L^2(\R^d)},
	\qquad
	d\geq 3
	,
\end{equation}
or equivalently
\begin{equation}\label{eq:KYequiv2}
	\big\||x|^{-1}f\big\|_{L^2(\R^d)}\leq C\||x|\Delta f\|_{L^2(\R^d)},
	\qquad
	d\geq3
	,
\end{equation}
which is a weighted version of the well-known Rellich inequality (see \cite{R}), namely
$$
\left\||x|^{-2}f\right\|_{L^2(\R^d)}\leq C\|\Delta f\|_{L^2(\R^d)},
\qquad
d\geq5.
$$
The proof by Kato and Yajima relies on the Fourier representation of solutions to \eqref{eq:euclidean}; nevertheless, Burq, Planchon, Stalker, and Tahvildar-Zadeh proved \eqref{eq:KY} in \cite{BPST1, BPST2} by the \textit{multipliers method}, which only involves algebraic manipulations of equation \eqref{eq:euclidean}. This also permits to involve in \cite{BPST1, BPST2} some repulsive potentials with critical behavior with respect to the free scaling in the topic.
Later, Barcel\'o, Vega, and Zubeldia proved in \cite{BVZ} a stronger estimate than \eqref{eq:KY}. More precisely, if $u$ is a solution to \eqref{eq:euclidean}, denoting by $\lambda=\lambda_1+i\lambda_2$, and
\begin{equation}\label{eq:upm_euc}
	u^\pm (x)
	:=
	e^{\pm i \sgn(\lambda_2) \sqrt{|\lambda_1|} |x|} u(x),
\end{equation}
where $\sgn$ is the sign function defined in \eqref{eq:defs}, the following estimate holds in the cone $|\lambda_2|\leq\lambda_1$:
\begin{equation}\label{eq:KR}
\|\nabla u^-\|_{L^2(\R^d)}\leq C\||x|f\|_{L^2(\R^d)},
\qquad
d\geq3.
\end{equation}
The zero-energy case has to be understood as a weighted Hardy--Rellich inequality (see for instance \cite{CCF}) for a suitable gauge transformation on $u$.
Notice that \eqref{eq:KR} scales as \eqref{eq:KY}, namely if $u_k(x):=u(\frac{x}{k})$ for $k\in\R$, then
$$
\left\||x|^{-1}u_k\right\|_{L^2(\R^d)}^2=k^{d-2}\left\||x|^{-1}u\right\|_{L^2(\R^d)}^2,
\qquad
\left\|\nabla u_k^-\right\|_{L^2(\R^d)}^2=k^{d-2}\left\|\nabla u^-\right\|_{L^2(\R^d)}^2.
$$
On the other hand, \eqref{eq:KR} implies \eqref{eq:KY} in dimension $d\geq3$, thanks to the Hardy inequality
\begin{equation}\label{eq:hardy}
\int_{\R^d}\frac{|u|^2}{|x|^2}
=
\int_{\R^d}\frac{|u^-|^2}{|x|^2}\leq \left( \frac{2}{d-2} \right)^2 \int_{\R^d}|\nabla u^-|^2,
\qquad
d\geq3.
\end{equation}

In dimensions $d=1,2$ estimate \eqref{eq:KR} fails, as well as \eqref{eq:KY}. Indeed, their validity would imply absence of eigenvalues for $-\Delta+V$ for sufficiently small $V$, while it is known (see \cite[Theorem~XIII.11]{ReedSimonIV}) that even compactly supported potentials $V$ in general produce eigenvalues (the Laplacian is said to be \textit{critical}). Although \eqref{eq:KR} can be understood from the Fourier analytical viewpoint as a radiation condition for solutions to the Helmholtz equation, the proof in \cite{BVZ} is inspired by the method introduced by Ikebe and Saito in \cite{IS} and completely relies on real analytical techniques based on multipliers methods. This permits the authors to involve electromagnetic perturbations in their statements, and obtain quite strong informations from the point of view of scattering theory. To complete the picture, in \cite{FKV} the authors used techniques similar to the ones in \cite{BVZ} to get spectral information for non-self-adjoint Hamiltonians.


\subsection{Aim of the paper}

Moving into the setting of the Heisenberg group, it is natural to investigate about the validity of inequalities in the same style as in \eqref{eq:KY}, \eqref{eq:KYequiv}, \eqref{eq:KYequiv2}, and \eqref{eq:KR}. As argued above, the Fourier analysis does not seem to help in understanding an estimate like \eqref{eq:KR}, and finding an analogous result in this way can be challenging. Therefore, we will here rely on the multipliers method.

Similarly to what happens for the Hardy inequality in the Heisenberg group, a first difficulty is to understand what weights should be involved in such estimates.
In fact, one may firstly wonder whether \eqref{eq:KYequiv2} is valid for {\it horizontal} weights, i.e. functions depending only on the horizontal variable $z$. Comparing with the horizontal Hardy inequality \eqref{eq:hardyDa} below, one may wonder if an inequality like
\begin{equation}\label{eq:KYH}
	\left\||z|^{-1}f\right\|_{L^2_{(z,t)}(\H^d)}\leq C\||z| \L f\|_{L^2_{(z,t)}(\H^d)}
\end{equation}
can hold for some positive constant $C$. If \eqref{eq:KYH} was true, then it would imply, by standard Kato-smoothing, the following estimate
\begin{equation}\label{eq:ksmo}
	\left\||z|^{-1}e^{-is\L}f\right\|_{L^2_sL^2_{(z,t)}(\H^d)}\leq C\|f\|_{L^2_{(z,t)}(\H^d)}
\end{equation}
for the Schr\"odinger evolution equation associated to $\L$.
Nevertheless, \eqref{eq:ksmo} cannot hold, since it would contradict the existence of soliton-like solutions to the evolution flow generated by $e^{-is\L}$, as showed by Bahouri, G\'erard and Xu in \cite{BGX} (see also \cite{BBG}, where some averaged Strichartz estimates are proved). We conclude that \eqref{eq:KYH} is necessarily false.

Another difficulty arises form the multipliers method itself: even though it has proved to be a robust technique valuable in providing spectral information in various contexts (see e.g. \cite{CFK,CFK23,CK,FKV}), unfortunately it lacks of an abstract formulation and it turns out to be difficult to generalize. 
Formally, the method relies on some precise algebraic manipulations of identities, obtained by the Helmholtz equation choosing suitable test functions (the multipliers of the method).
As a consequence, it is a quite delicate machinery and usually strongly sensitive to the model under consideration.

In our main Theorem\til\ref{thm:grad-est}, we succeed in finding the analogous of the estimates \eqref{eq:KY} and  \eqref{eq:KR} in the context of the Heisenberg group. However, due to the aforementioned technical difficulties arising from the multipliers method, only under radial assumptions 
for the weak solutions of \eqref{eq:prob0}.
Leveraging on this result, in Theorem\til\ref{thm:V1} we are able to exclude the existence of radial eigenfunctions for suitable small, possibly complex, perturbation of the sublaplacian.
Theorems\til\ref{thm:pp},\til\ref{thm:grad-est2} and\til\ref{thm:V2} are dedicated instead to improve the first two results in the presence of real and/or repulsive potentials.
It is not clear if a different choice of the multipliers, or of the change of gauge $u^-$, could let us approach also the general case, or if on the contrary the multipliers method is not the most suitable one for the Heisenberg group. 
Here we start a research program aiming to produce resolvent estimates, particularly weighted $L^2$-estimates, in $\H^d$.
Different approaches will be explored by the authors in a subsequent work \cite{FMRS}. 

This paper is structured as follows. In the next Section\til\ref{sec:main_results}, after recalling some basic knowledge on the Heisenberg group and the sublaplacian, we state our main results. Section\til\ref{sec:Hardy} is devoted in recalling (weighted) Hardy inequalities, crucial tools in our demonstration. Section\til\ref{sec:proof} is dedicated to the proof of Theorem\til\ref{thm:grad-est}, which constitutes a blueprint for the proofs of Theorems\til\ref{thm:V1}--\ref{thm:V2}, considered in Sections\til\ref{sec:proof2}--\ref{sec:proof4}.


\section{Main results}\label{sec:main_results}

In dimension $d\geq1$, let $\H^d:=\C^d\times\R$ be the $d$-dimensional Heisenberg group
endowed with the usual metric structure induced by the Koranyi norm
\begin{equation}
	\label{eq:koranyi}
	|(z,t)|_{\H} := (|z|^4 + t^2)^{1/4},
	\qquad
	(z,t)\in\H^d,
\end{equation}
as well as the group law
\begin{equation*}
	(z,t)(z',t'):=\big(z+z',t+t'+2\,\im(z\cdot \overline{z'})\big),
\end{equation*}
being $z\cdot \overline{z'}=z_1\overline{z_1'}+\cdots+z_d\overline{z'_d}$ with $z,z'\in\C^d$. We will regard $\H^d$ as a measure space together with the Haar measure, i.e. nothing else than the Lebesgue measure $dz\,dt$ on $\C^d\times\R$. 
From now on we will make use of the natural identification $\C^d\simeq\R^{2d}$, thanks to which we can write
$$
\C^d\ni (z_1,\dots,z_d) = z \simeq(x,y)=(x_1,\dots,x_d,y_1,\dots,y_d) \in \R^{2d},
$$
where
$$
x_j=\re z_j,\quad y_j=\im z_j,
\quad
|z|^2=|x|^2+|y|^2.
$$
Hence we can introduce, for $j \in \{1,\dots,d\}$, the $(2d+1)$ left-invariant vector fields 
\begin{equation*}
	X_j := \frac{\partial}{\partial x_j} + 2{y_j} \frac{\partial}{\partial t},
	\qquad
	Y_j := \frac{\partial}{\partial y_j} - 2{x_j} \frac{\partial}{\partial t},
	\qquad
	T := \frac{\partial}{\partial t}
	,
\end{equation*}
and the span $\mathcal D:=\text{span}\{X_j,Y_j\}_{j=1,\dots,d}$ which provides a sub-Riemannian structure on $\R^{2d+1}$. 
The above vector fields satisfy, for any $j,k = 1, \dots, d$, the commutation relations
\begin{equation}\label{eq:commutators}
	[X_j,X_k] = [Y_j,Y_k] = 0, 
	\qquad
	[X_j,Y_k] = - 4 \delta_{jk} T,
\end{equation}
where $\delta_{jk}$ is the Kronecker's delta.
The associated Laplacian, called
{\it sublaplacian}, is the hypoelliptic operator $\L$  defined by
\begin{equation*}
	\L := -\sum_{j=1}^{d} \big(X_j^2 + Y_j^2\big).
\end{equation*}
Using the notation
\begin{equation*}
	\grad := (X_1,\dots,X_d,Y_1,\dots,Y_d),
\end{equation*}
the sublaplacian can be written in the divergence form $\L = - \grad\cdot\grad$. 
We will refer to the vector $\grad$ as the \textit{horizontal gradient}, since it does not involve the vertical direction generated by $T=\frac{\partial}{\partial t}$, whereas $\nabla=(\grad,T)$ is the (complete) \textit{gradient} in the Heisenberg group. 

Let us recall also the natural dilations $\delta_\lambda$ on the Heisenberg group, defined by
\begin{equation}\label{def:dilation}
	\delta_\lambda(z,t) = (\lambda z , \lambda^2 t).
\end{equation}
Notice that this dilation has a good interaction with the sublaplacian and the Haar measure, in the sense that, taking (for instance) a smooth compactly supported function $f\in \mathcal{C}_0^{\infty}(\H^d)$, then
\begin{equation*}
	\L(f \circ \delta_\lambda) = \lambda^2 (\L f) \circ \delta_\lambda.
\end{equation*}
This means that $\L$ is \emph{homogeneous of degree 2}, and for the Haar measure $\text{d}$ we have
\begin{equation*}
	\text{d}\,(\delta_\lambda(z,t)) = \lambda^Q dzdt
\end{equation*}
where $Q=2d+2$ is the \emph{homogeneous dimension} of $\H^d$.

We also introduce the Sobolev space $H^1(\H^d)$ of the functions $u:\H^d\to\C$ such that
\begin{equation*}
	\norm{u}_{H^1(\H^d)}
	:=
	\left( \norm{u}^2+\norm{\grad u}^2 \right)^{1/2}
	<
	+\infty
	%
\end{equation*}
where $\norm{\cdot}$ denotes the $L^2$-norm, viz.
\begin{equation*}
	\norm{u} 
	\equiv 
	\norm{u}_{L^2(\H^d)}
	:=
	\left(
	\int_{\H^d}
	|u|^2
	\,dz\,dt
	\right)^{1/2}
	.
\end{equation*}
Following the nomenclature, e.g. in \cite{GL}, in the present work a function $g \colon \H^d \to \C$ is defined to be \textit{radial}\footnote{Notice that the literature is not unanimous in the use of this terminology: sometimes a function $g \colon \H^d \to \C$ is said to be radial if there exists $g_0 \colon [0,+\infty) \times \R \to \C$ such that $g(z,t)=g_0(|z|,t)$.} if there exists some function $g_0 \colon [0,+\infty) \to \C$ such that
\begin{equation}\label{eq:radial_def}
	g(z,t) = g_0(\Kor).
\end{equation}
For a complete description of $\H^d$, we address the reader to the standard reference \cite{F2}.


Now, given $\lambda=\lambda_1+i\lambda_2\in\C$, and $u \colon \H^d\to\C$, we define
\begin{equation}\label{eq:upm}
	u^\pm (z,t)
	:=
	e^{\pm i \psi(z,t)} u(z,t),
	\qquad
	\psi(z,t)
	:=
	\sqrt{\frac{d}{2}} \frac{\Gamma\left(\frac{d}{2}\right)}{\Gamma\left(\frac{d+1}{2}\right)} \sgn(\lambda_2) \sqrt{|\lambda_1|} 
	\Kor
	,
\end{equation}
where $\Gamma$ is the Euler Gamma function and
\begin{equation}\label{eq:defs}
	\sgn(w):=
	\begin{cases}
		\frac{w}{|w|} &\text{if $w\neq0$,}
		\\
		1 &\text{otherwise}.
	\end{cases}
\end{equation}
We are finally ready to state the first result of this paper.

\begin{theorem}\label{thm:grad-est}
	Let $d\ge1$ and let $u\in H^1(\H^d)$ be a radial weak solution of \eqref{eq:prob0}. Then, for any $\delta>0$, the following estimates hold:
	\begin{equation}\label{eq:est1}
		\norm{\nabla_{\!\H} u} 
		\le \frac{1+1/\delta}{d} \norm{ \frac{\Kor^2}{|z|} f}, \qquad
		\text{if $|\lambda_2|>\delta \lambda_1 $}
	\end{equation}
	and
	\begin{equation}\label{eq:est2}
		\norm{\nabla_{\!\H}u^-}
		\le
		K_{d}(\delta)
		\norm{ \frac{\Kor^2}{|z|} f }, \qquad
		\text{if $|\lambda_2|\leq \delta \lambda_1 $}
	\end{equation}
	where
	\begin{equation}\label{eq:kd}
		K_{d}(\delta)
		:=
		\min_{\gamma>0} 
		\left(
		\frac{8d+2+\gamma\sqrt{\delta}}{4d}+\sqrt{\frac{(8d+2+\gamma\sqrt{\delta})^2}{16d^2}+\frac{\sqrt{\delta}}{2\gamma}} 
		\right)
		.
	\end{equation}
	In particular, the uniform estimate 
	\begin{equation}\label{eq:katoyajima}
		\norm{ \frac{|z|}{\Kor^2} u} \le  \kappa_d \norm{ \frac{\Kor^2}{|z|} f }
	\end{equation}
	holds, where 
	\begin{equation}
		\label{kappa}
		\kappa_d := \min_{\delta>0} \max\left\{ \frac{1+1/\delta}{d^2}, \frac{K_{d}(\delta)}{d} \right\}.
	\end{equation}
\end{theorem}

\begin{remark}\label{rem:000}
	It is interesting to compare these estimates to their Euclidean counterparts, and to the well-known Hardy-type inequality by Garofalo and Lanconelli in \cite[Corollary 2.1]{GL}, on which we will return later in Section \ref{sec:Hardy}.
	First of all, notice that the weight $\frac{\Kor^2}{|z|}$ is equal to the reciprocal of the weight appearing in the Garofalo--Lanconelli inequality. The geometrical meaning of this quantity is soon explained: it is the ratio between the Koranyi gauge and the modulus of its horizontal gradient (see \eqref{eq:kder1} and \eqref{eq:kder2}). Explicitly:
	\begin{equation*}
		\frac{\Kor^2}{|z|} = \Kor \cdot \left( \frac{|z|}{\Kor} \right)^{-1} = \frac{\Kor}{|\grad\Kor|}
		.
	\end{equation*}
	In the Euclidean case, replacing $\Kor$ with $|x|$ and $\grad$ with $\nabla$ in the right-hand side, this quantity simply reduces to the weight $|x|$ appearing on the right-hand side of \eqref{eq:KY}.
	Formally, comparing \eqref{eq:KY} and \eqref{eq:katoyajima} (but also the Hardy inequalities \eqref{eq:hardy} and \eqref{eq:hardyGL}), one can switch from the Euclidean to the Heisenberg case of these inequalities according to the table:
	\begin{table}[H]
		\centering
		\label{tab0}
		\begin{tabular}{c@{}c@{}c}
			\toprule
			Euclidean & \qquad\quad\quad & Heisenberg\\ \midrule
			$d$ & & $Q=2d+2$ \\
			$|x|$ & & $\Kor$ \\
			$\nabla$ & & $\grad$ \\
			$|\nabla|x||=1$ & & $|\grad\Kor|=\frac{|z|}{\Kor}$ \\
			\bottomrule
		\end{tabular}
	\end{table}
	Clearly this comparison, even if evocative, must be taken with (more than) a grain of salt, since the non-flat geometry of the Heisenberg group and its intrinsic structure make the situation much more peculiar. For example, and this is crucial in our proof, whereas the Euclidean radial derivative can be written as a directional derivative, $\partial_\rho = \frac{\nabla|x|}{|\nabla|x||} \cdot \nabla$, this is not true in the Heisenberg group, where the natural radial derivative is related to the Euler vector field defined in \eqref{eq:EulerVF} below. This motivates the introduction of a change of coordinates in the Heisenberg group (as we will do later, see Subsection \ref{sub:change}). 
	Finally, it should be noted that our result holds for any $d\ge1$, whereas the Euclidean one holds for $d\ge3$.
	These dimensional conditions are inherited by the respective ones on the Hardy's inequalities (the Euclidean Laplacian is critical for $d=1,2$, and the Euclidean Hardy's inequality is not valid in these dimensions).
	Again, one can see that the Euclidean condition $d\ge3$ becomes $Q\ge3$ (and so $d\ge1$) in the Heisenberg case.
\end{remark}

\begin{remark}\label{rem:100}
The Euclidean version of Theorem \ref{thm:grad-est} is contained inside the proof of \cite[Theorem 8]{FKV}, although it is not stated as an a priori estimate. In addition, no optimization on the opening parameter $\delta$ of the cone  $|\lambda_2|\le \delta |\lambda_1|$ is provided in \cite{FKV}, where merely $\delta=1$. We highlight that here we are giving explicit constants in the estimates \eqref{eq:est1}, \eqref{eq:est2}, and \eqref{eq:katoyajima}, unlike \cite{BVZ, BPST1, BPST2,FKV}. 
This does not depend on our radial assumption on the solution, and these constants are valid also in the Euclidean case, for any $u \in \H^1(\R^d)$, after replacing $d$ with $d/2-1$ in the above formulas (see next Remarks\til\ref{rem:200} and\til\ref{rem:500} too).
Up to our knowledge, this kind of result was totally unknown in the setting of the Heisenberg group.
\end{remark}

\begin{remark}\label{rem:200}
By elementary calculus, it is easy to check that the minimum in the expression of $K_{d}(\delta)$ is attained when $\gamma=\gamma_\delta$, where $\gamma_\delta>0$ is the unique positive solution of the cubic equation 
\begin{equation}
\label{eq:cubica}
	\sqrt{\delta} \gamma_\delta^3 + (4d+1) \gamma_\delta^2 - d^2 = 0.
\end{equation}
Therefore, we have
\begin{equation*}
	K_{d}(\delta) = \frac{d}{\gamma_\delta^2},
\end{equation*}
or equivalently, $K_d(\delta)$ is the unique positive solution of the equation
\begin{equation}
	\label{eq:Kd-impdef}
	\sqrt{d \delta } K_{d}(\delta)^{-3/2} + (4d+1) K_{d}(\delta)^{-1} = d.
\end{equation}

As for $\kappa_d$, since $\frac{1+1/\delta}{d^2}$ is strictly decreasing in $\delta$, while $\frac{K_{d}(\delta)}{d}$ is strictly increasing with respect to $\delta$, then there exists a unique positive $\delta^*$ such that 
\begin{equation}
	\label{eq:kd-def2}
	\k_d = \frac{1+1/\delta_*}{d^2} = \frac{K_{d}(\delta_*)}{d}.
\end{equation}
From \eqref{eq:cubica} and \eqref{eq:Kd-impdef} we see that $(\gamma_{\delta_*},\delta_*)$ is the unique positive solution of the system
\begin{equation*}
	\begin{cases}
		\frac{1+1/\delta}{d^2} = \frac{1}{\gamma^2},
		\\
		\sqrt{\delta} \gamma^3 + (4d+1) \gamma^2 - d^2 = 0
	\end{cases}
\end{equation*}
which leads to 
\begin{gather}
	\label{eq:delta*}
	\frac{\delta_*^2}{\sqrt{1+\delta_*}} + 4\delta_* = \frac{1}{d}, \qquad \gamma_{\delta_*}=\frac{d \delta_*}{\delta_*+1}.
\end{gather}
Equivalently, $\kappa_d$ is the unique positive solution of the equation
\begin{equation}
\label{kappanueva}
	\frac{\kappa_d^{-2}}{\sqrt{d^2- \kappa_d^{-1}}} + (4d+1) \kappa_d^{-1} = d^2.
\end{equation}
In the following table we show the explicit value of the constants $\delta_*$ and $\kappa_d$ in some particular cases (the values are approximated to the fifth decimal place).
\begin{table}[H]
	\centering
	\label{tab1}
	\begin{tabular}{c@{}ccc@{}c}
		\toprule
		$d$ & \,\, & $\delta_*$               & \,\, & $\kappa_d$ \\ \midrule
		$1$ &  &  $2.37340 \cdot 10^{-1}$ &  & $5.21337$ \\
		$2$ &  &  $1.21514 \cdot 10^{-1}$ &  & $2.30737$ \\
		$3$ &  &  $8.17278 \cdot 10^{-2}$ &  & $1.47064$ \\
		$4$ &  &  $6.15799 \cdot 10^{-2}$ &  & $1.07744$ \\
		$5$ &  &  $4.94043 \cdot 10^{-2}$ &  & $8.49645 \cdot 10^{-1}$ \\
		\bottomrule
	\end{tabular}
\end{table}
Again we stress that an explicit value for $\kappa_d$ was not known even in the Euclidean case. To obtain the Euclidean analogue it is enough to replace $d$ with $d/2-1$ above. It is an interesting open question to establish whether the constant $\kappa_d$ is optimal. Finally, it can be observed that $\delta_*$ and $\kappa_d$ are strictly decreasing to $0$ as $d \to +\infty$, as per their implicit definition. Also, $\k_d$ satisfies the na\"{\i}f bound from below
\begin{equation}
	\label{eq:kd-bound}
	\kappa_d > \frac{4d+1}{d^2} 
\end{equation} 
since from its implicit definition \eqref{kappanueva} one has $d^2-(4d+1) \k_d^{-1} >0$.
\end{remark}


As a consequence of Theorem \ref{thm:grad-est}, we can prove that the eigenfunctions of some complex perturbations of $\L$ can not be radial, in the sense that if $u \in H^1(\H^d)$ is a radial solution of $(\L+V)u=\lambda u$ for some $\lambda\in\C$, then $u\equiv0$. 

Firstly, we recall that $\L$ is the positive self-adjoint operator on $L^2(\H^d)$ associated to the quadratic~form
$$
q_0[\psi]:=\int_{\H^d}|\nabla_{\!\H}\psi|^2,
\qquad
D(q_0)=H^1(\H^d).
$$
The spectrum of $\L$ is purely continuous and given by
$
\sigma(\L)=\sigma_c(\L)=[0,+\infty) ,
$
see e.g. \cite{DMW} and the references therein. Let now $V:\H^d\to\C$ be a measurable function and assume $V$ is subordinated to $\L$, with bound less than one, namely
\begin{equation}\label{eq:subordination}
\exists \, a<1 \quad \text{ such that }
\qquad
\int_{\H^d}|V||\psi|^2\leq a\int_{\H^d}|\nabla_{\!\H}\psi|^2,
\end{equation}
for all $\psi\in H^1(\H^d)$. This implies that the quadratic form
$$
q_V[\psi]:=\int_{\H^d}V|\psi|^2,
\qquad
D(q_V):=\Big\{\psi\in L^2(\H^d)\colon \int_{\H^d}|V||\psi|^2<\infty\Big\}
$$
is relatively bounded with respect to $q_0$, with bound less than 1. Consequently, one can define an \mbox{$m$-sectorial} differential operator $\L+V$ associated to the quadratic form $q:=q_0+q_V$, since this is a closed form (see \cite[Chapter VI, Theorem 2.1]{K}). The subordination condition \eqref{eq:subordination} needs to be understood in terms of the Hardy inequalities on $\H^d$, to which we devote Section~\ref{sec:Hardy} below.

We have the following consequence of Theorem\til\ref{thm:grad-est}. 

\begin{theorem}\label{thm:V1}
	Let $d\geq1$, let $V\colon \H^d\to\C$ be a measurable function satisfying \eqref{eq:subordination}.
	Assume that
	\begin{equation}\label{eq:assV1}
		\int_{\H^d} \frac{|(z,t)|_{\H}^4}{|z|^2} |V|^2|\psi|^2\leq b^2\int_{\H^d}|\nabla_{\!\H}\psi|^2,
	\end{equation}
	for all $\psi\in H^1(\H^d)$ and for some positive constant $b>0$ satisfying
	\begin{equation}\label{eq:nuovanuova}
		b< \frac{1}{d \, \kappa_d}
	\end{equation}
	where $\kappa_d$ is given by \eqref{kappa} or equivalently \eqref{kappanueva}.
	Then, $\L+V$ has no radial eigenfunctions.
\end{theorem}

\begin{remark}\label{rem:new}
When $u=u_0(|(z,t)|_\H)$ is a radial function, from an easy computation we get
$$
-\mathcal Lu(z,t)=\frac{|z|^2}{|(z,t)|^2_\H}\left(u''_0(|(z,t)|_\H)+\frac{2d+1}{|(z,t)|_\H }u'_0(|(z,t)|_\H)\right),
$$
showing that the sublaplacian does not preserve radial symmetry. On the one hand, this is a weakness of  Theorem \ref{thm:V1}, since it shows that there is no reason for which radial eigenfunctions should exist, even for radial potentials $V$. On the other hand, this opens an interesting question about the possible symmetries of the eigenfunctions of $-\mathcal L+V$, which will be object of further investigation in the next future, with a special attention devoted to the case of cylindrically symmetric potentials.
\end{remark}
\begin{remark}\label{rem:500}
Comparing the above result with the Euclidean analogue in \cite[Theorem 2]{FKV} for the Schr\"odinger Hamiltonian, we see that the bound \eqref{eq:nuovanuova} is better. Indeed, if one replaces $d$ by $d/2-1$ then the above condition reads as $b < \frac{2}{d-2} \, \kappa_{d/2-1}^{-1}$, which is better than the bound $b < \frac{d-2}{5d-8}$ in \cite[Theorem 2]{FKV}, as it easily follows from analogous arguments as before.

\end{remark}

\begin{remark}\label{rem:550}
By Cauchy--Schwarz and the Hardy inequality \eqref{eq:hardyGL}, one can estimate
$$
\int_{\H^d}|V||\psi|^2\leq
\left(\int_{\H^d} \frac{\Kor^4}{|z|^2} |V|^2|\psi|^2\right)^{\frac12}
\left(\int_{\H^d} \frac{|z|^2}{\Kor^2} \cdot \frac{|\psi|^2}{\Kor^2}\right)^{\frac12}
\leq\frac b{d}\int_{\H^d}|\nabla_{\!\H}\psi|^2,
$$
and clearly \eqref{eq:subordination} is automatically satisfied, if we assume  \eqref{eq:nuovanuova} and we take  \eqref{eq:kd-bound} into account.
To obtain explicit examples of potentials $V$ satisfying \eqref{eq:assV1}, it is sufficient to involve the Hardy inequality \eqref{eq:hardyGL} to realize that the critical potential scaling is given by homogeneous functions like $|(z,t)|_{\H}^{-2}$.
\end{remark}

We now pass to some generalizations of Theorems \ref{thm:grad-est} and \ref{thm:V1}, in which some symmetric perturbation of order zero is also present, in order to detect a suitable smallness condition which fits with the physical notion of {\it repulsivity}. 
	Before doing this, we need to introduce the Euler vector field defined by
	\begin{equation}\label{eq:EulerVF}
		\mathbb{E} := z \cdot \nabla_z + 2 t \partial_t,
	\end{equation}
	which characterizes the homogeneity on the Heisenberg group, in the sense that it is the only vector field such that $\mathbb{E} f = \nu f$ if $f$ is $\nu$-homogeneous with respect to the Heisenberg dilation \eqref{def:dilation}, explicitly $f \circ \delta_\lambda = \lambda^\nu f$. Then, we define the operator
	\begin{equation*}
		\partial_\kor := \Kor^{-1} \mathbb{E}
	\end{equation*} 
	which arises naturally in the literature as the radial derivative, in the Heisenberg setting, taken with respect to the Koranyi gauge, and hence as a counterpart of the radial derivative $\frac{x}{|x|} \cdot \nabla$ from the Euclidean case (see for instance \cite[Section\til3.1]{RS}, but also Remark\til\ref{rmk:spherical} below).
Let us consider now the equation
\begin{equation}\label{eq:prob2}
-\L u-Vu+\lambda u=f,
\end{equation}
where $u,f\colon \H^d\to\C$, $\lambda\in\C$, and the potential $V\colon\H^d\to\R$ is a real-valued measurable function on $\H^d$. The weak formulation of \eqref{eq:prob2} is now
\begin{equation}\label{eq:solution2}
	- 
	\langle \grad v, \grad u \rangle
	-
	\langle v,Vu\rangle
	+
	\lambda
	\langle v, u \rangle
	=
	\langle v, f \rangle.
\end{equation}
Our next result is a generalization of Theorem \ref{thm:grad-est}, related to equation \eqref{eq:prob2}. We start with the case of a positive potential. Below, $W^{k,d}_{\operatorname{loc}}(\H^d;\R)$ stands for the usual Sobolev space of real-valued functions $f\in L^d(\H^d)$ such that $f$ and its weak derivatives up to order $k$ have a finite $L^d$-norm.

\begin{theorem}[Positive potentials]\label{thm:pp}
	Let $d\ge1$, $u\in H^1(\H^d)$ be a radial weak solution of \eqref{eq:prob2}, with non-negative potential $V\in L^1_{\operatorname{loc}}(\H^d;\R)\cap W^{1,d}_{\operatorname{loc}}(\H^d;\R)$, and
	assume that there exists a constant
	$0 \le b < 1$ such that, for any $\psi\in H^1(\H^d)$, we have
	\begin{equation}\label{eq:b2b}
		\int_{\H^d} \left[ \partial_{\kor}(\Kor V) \right]_+|\psi|^2  \leq b^2\int_{\H^d}|\nabla_{\!\H}\psi|^2,
	\end{equation}
	where $[\cdot]_+=\max\{0,\cdot\}$.
	Then, for any $\delta>0$, the following estimates hold:
	\begin{equation}\label{eq:est3b}
		\norm{\grad u} 
		\le \frac{1+1/\delta}{d} \norm{ \frac{\Kor^2}{|z|} f}
		\qquad
		\text{if $|\lambda_2|>\delta\lambda_1$,}
	\end{equation}
	and
	\begin{equation}\label{eq:est4b}
		\norm{\grad u^-}
		\le
		K_{d,b}(\delta)
		\norm{ \frac{\Kor^2}{|z|} f} 
		\qquad
		\text{if $|\lambda_2|\leq\delta\lambda_1$,}
	\end{equation}
	where
	\begin{equation}\label{eq:kdb}
		K_{d,b}(\delta)
		:=
		\min_{\gamma>0} 
		\left(
		\frac{8d+2+\gamma\sqrt{\delta}}{4d(1-b^2)}
		+
		\sqrt{ \left[\frac{8d+2+\gamma\sqrt{\delta}}{4d(1-b^2)}\right]^2 +\frac{\sqrt{\delta}}{2\gamma (1-b^2)}}
		\right)
		.
	\end{equation}
	
	In particular, the following uniform estimate holds:
	\begin{equation}\label{eq:katoyajima2}
		\norm{ \frac{|z|}{\Kor^2} u} \le \kappa_{d,b} \norm{ \frac{\Kor^2}{|z|} f}
	\end{equation}
	where 
	\begin{equation}\label{eq:md2}
		\kappa_{d,b} :=
		\min_{\delta>0}
		\max\left\{\frac{1+1/\delta}{d^2}\,,\,\frac{K_{d,b}(\delta)}{d}\right\},
	\end{equation}
	and $K_{d,b}(\delta)$ is the constant in \eqref{eq:kdb}.
\end{theorem} 

\begin{remark}[Repulsive potentials]
Notice that if $V\ge0$ is \textit{repulsive}, in the sense that $\partial_\kor(\Kor V)\le0$, then $b=0$ in \eqref{eq:b2b} and the above result reduces to Theorem\til\ref{thm:grad-est}.
Notice also that $\partial_\kor(\Kor V) = V + \Kor \partial_\kor V = (1+\mathbb{E}) V$.
\end{remark}

In the following result we allow a negative part for $V$.
\begin{theorem}\label{thm:grad-est2}
	Let $d\ge1$, $u\in H^1(\H^d)$ be a radial weak solution of \eqref{eq:prob2}, with potential $V\in L^1_{\operatorname{loc}}(\H^d;\R)\cap W^{1,d}_{\operatorname{loc}}(\H^d;\R)$, and
	assume that there exist two constants
	$0 \le b_1,b_2 < 1$ such that, for any $\psi\in H^1(\H^d)$, we have
	\begin{equation}
		\int_{\H^d}V_-|\psi|^2  \leq b_1^2\int_{\H^d}|\nabla_{\!\H}\psi|^2
		\label{eq:b1}
	\end{equation}
	and
	\begin{equation}
		\int_{\H^d} [\partial_\kor(\Kor V)]_+|\psi|^2  \leq b_2^2\int_{\H^d}|\nabla_{\!\H}\psi|^2,
		\label{eq:b2}
	\end{equation}
	where $V_-=\max\{0,-V\}$.
	Then, for any $\delta>0$, the following estimates hold:
	\begin{equation}
		\norm{\grad u} 
		\le \frac{1+1/\delta}{d(1-b_1^2)} \norm{ \frac{\Kor^2}{|z|} f}
		\qquad
		\text{if $|\lambda_2|>\delta\lambda_1$,}
		\label{eq:est3}
	\end{equation}
	and
	\begin{equation}
		\norm{\grad u^-}
		\le
		M_{d,b_2}(\delta)
		\norm{ \frac{\Kor^2}{|z|} f}
		\qquad
		\text{if $|\lambda_2|\leq\delta\lambda_1$,}
		\label{eq:est4}
	\end{equation}
	where
	\begin{equation}\label{eq:kd2}
		M_{d,b_2}(\delta)
		:=
		\min_{\substack{\gamma_1>0 \\ 0<\gamma_2<1}} 
		g_{d,\delta,b_2}(\gamma_1,\gamma_2)
	\end{equation}
	and 
	$g_{d,\delta,b_2}(\gamma_1,\gamma_2):=a+\sqrt{a^2+b}$, with $$a=\frac{ 
			4d+1 + \frac{\gamma_1 \sqrt{\delta}}{2}
			+
			\frac{1}{1-b_2^2} \left(\frac{\sqrt{\delta}}{8d \gamma_2} \right)^2 
		}
		{2d(1-b_2^2)(1-\gamma_2^2)},\qquad b=\frac{\sqrt{\delta}}{2\gamma_1(1-b_2^2)(1-\gamma_2^2)}.
$$	
In particular, the following uniform estimate holds:
	\begin{equation}\label{eq:katoyajima2-2}
		\norm{ \frac{|z|}{\Kor^2} u} \le \mu_{d,b_1,b_2} \norm{ \frac{\Kor^2}{|z|} f}
	\end{equation}
	where 
	\begin{equation}\label{eq:md2-2}
		\mu_{d,b_1,b_2} :=
		\min_{\delta>0}
		\max\left\{\frac{1+1/\delta}{d^2(1-b_1^2)}\,,\,\frac{M_{d,b_2}(\delta)}{d}\right\},
	\end{equation}
	and $M_{d,b_2}(\delta)$ is the constant in \eqref{eq:kd2}.
\end{theorem} 

\begin{remark}
Observe that if $b_1=b_2=0$ in  Theorem \ref{thm:grad-est2}, we are not recovering the results in the free case. This is due to the presence of an additional term involving the negative part of $V$ which scales differently than the other terms, see \eqref{eq:idue0} below. 
The corresponding  minimization problem of the constants is quite more involved, see Remark\til\ref{rem:1000} below. 

\end{remark}

\begin{remark}\label{rem:800}
The regularity assumption $V\in L^1_{\text{loc}}(\H^d)$ is the minimal requirement in order to give a meaning to the weak formulation in \eqref{eq:solution2}. On the other hand, the assumption $V\in W_{\operatorname{loc}}^{1,d}(\H^d)$ is needed in order to produce a rigorous proof, after a regularization argument. We send the reader to the proof of \cite[Theorem 3.1]{CFK}, in which this argument is treated in detail, in the Euclidean setting (which does not present relevant differences from our case, from the point of view of regularity).
\end{remark}
\begin{remark}\label{rem:900}
Since $V$ is real, we only need to assume here the subordination condition on the negative part \eqref{eq:b1} on $V_-$, instead of \eqref{eq:subordination}. Indeed, the Friedrichs' extension theorem allows to define the self-adjoint differential operator $\L+V$ associated to the natural quadratic form, which under assumption \eqref{eq:b1} is positive.
\end{remark}

\begin{remark}\label{rem:1000}
The shape of the constant {$M_{d,b_2}(\delta)$} is more complicated, compared with the free constant $K_d$ in \eqref{eq:kd}, since it involves a double optimization problem. The difference arises since we need to estimate some terms involving $V$ (see \eqref{eq:idue0}) which live at a different scale with respect to the one given by estimate \eqref{eq:est4}. We stress that this problem does not exist if some additional sign condition is assumed (e.g., if $V\geq0$, so that the term to estimate in \eqref{eq:idue0} is null). In conclusion, we decided to write \eqref{eq:kd2} in such a general form, although the difficulty can be reduced in some specific cases. 
\end{remark}

We finally state our last theorem, which should be compared with Theorem\til\ref{thm:V1}, and where a radial condition on the real part of the potential is also present.

\begin{theorem}\label{thm:V2}
Let $d\geq1$, $V\in L^1_{\operatorname{loc}}(\H^d;\C)\cap W^{1,d}_{\operatorname{loc}}(\H^d;\C)$, and assume that there exist non-negative numbers $0\le b_1,b_2 <1$ and $b_3$
satisfying
\begin{equation}\label{eq:assf.Lambdanew}
	0 \le b_3
	<
	- \frac{1}{8d} - \frac{4d+1}{2} \sqrt{1-b_1^2}
	+
	\sqrt{ \left[ \frac{1}{8d} + \frac{4d+1}{2} \sqrt{1-b_1^2} \right]^2 + d (1-b_2^2) \sqrt{1-b_1^2} }
\end{equation}
such that, 
for all $\psi\in H^1(\H^d)$,
\begin{equation}
  \label{eq:assV12}
  \int_{\mathbb H^d} [\re V]_- \, |\psi|^2 
  \leq
  b_1^2\int_{\mathbb H^d}|\nabla_{\!\H}\psi|^2 \,,
\end{equation}
\begin{equation}
  \label{eq:asspart}
  \int_{\mathbb H^d} [\partial_\kor( \Kor \re V) ]_+ \, |\psi|^2
    \leq
   b_2^2\int_{\mathbb H^d}|\nabla_{\!\H}\psi|^2 \,,
   \end{equation}
   and
  \begin{equation}
   \int_{\mathbb H^d} \frac{\Kor^4}{|z|^2} \, |\im V|^2 \, |\psi|^2 
  \leq
  b_3^2\int_{\mathbb H^d}|\nabla_{\!\H}\psi|^2 \,
  .
  \label{eq:assV2}
  \end{equation}
 Then $\L+V$ has no radial eigenfunctions.
\end{theorem}

\begin{remark}\label{rem:2000}
Theorem \ref{thm:V2} is analogous to \cite[Theorem 3]{FKV} for the Euclidean case.
For the regularity requirement $V\in L^1_{\text{loc}}(\H^d;\C)\cap W^{1,d}_{\text{loc}}(\H^d;\C)$, the same discussion as in Remark \ref{rem:800} is valid.
We stress that one can ask the additional condition
\begin{equation*}
	b_1^2 + b_2^2 + \frac{b_3}{d} < 1
\end{equation*}
to ensure that $V$ is subordinated, and hence to have the operator $\L+V$ well-defined. Indeed, \eqref{eq:assV12} and \eqref{eq:asspart} in particular imply
	\begin{equation*}
		\int_{\H^d} |\re V| |\psi|^2 \le (b_1^2+b_2^2) \int_{\H^d} |\grad \psi|^2
	\end{equation*} 
whereas from \eqref{eq:assV2} and the Hardy inequality \eqref{eq:hardyDa} it follows that
	\begin{equation}
	\label{eq:useful}
		\int_{\H^d} |\im V| |\psi|^2 \le 
		\left( \int_{\mathbb H^d} \frac{\Kor^4}{|z|^2} \, |\im V|^2 \, |\psi|^2 \right)^{\frac12}
		\left(\int_{\H^d} \frac{|z|^2}{\Kor^4} {|\psi|^2} \right)^{\frac12}
		\le 
		\frac{b_3}{d} \int_{\H^d} |\grad \psi|^2.
	\end{equation}
\end{remark}

	\begin{remark}\label{rem:1100}
		It is relevant to stress that the kind of radiality assumption for the solutions in the theorems is a technical requirement introduced in order to apply the method of multipliers in the Heisenberg setting, as will be clear later from our proofs. In the same way, the definition of radial function as $g(z,t)=g_0(|(z,t)|_{\H})$ is more appropriate than the definition as $g(z,t)=g_0(|z|,t)$ for which, again, the methods and techniques in this paper do not work. The radial hypothesis can be relaxed if we stay outside the cone around the continuous spectrum of $\L$: \eqref{eq:est1}, \eqref{eq:est3b} and \eqref{eq:est3} holds even in the generic case, and perturbations of the sublaplacian, under the smallness assumptions in Theorem\til\ref{thm:V1} and\til\ref{thm:V2}, do not create eigenvalues far away from the non-negative half-line (since there are no non-vanishing solutions to $(\L+V)u=\lambda u$ for $|\lambda_2|>\delta \lambda_1$ with suitable $\delta>0$, see Sections\til\ref{sec:proof2} and\til\ref{sec:proof4}). 
		
		The real deal is the conic region $|\lambda_2| \le \delta \lambda_1$: do estimates \eqref{eq:est2}, \eqref{eq:est4b} and \eqref{eq:est4} hold even for non-radial solutions (possibly with a change of gauge different from \eqref{eq:upm})? 
		Or at least the Kato--Yajima-type estimates \eqref{eq:katoyajima}, \eqref{eq:katoyajima2} and \eqref{eq:katoyajima2-2}?
		The answer to these questions will be object of further investigation in a subsequent work \cite{FMRS}.
	\end{remark}

\section{Hardy inequalities on $\H^d$}\label{sec:Hardy}

We recall in this section some classic results on Hardy and weighted Hardy inequalities in the $\H^d$ setting, which will be useful in the sequel. First of all, we mention the celebrated work by Garofalo and Lanconelli \cite{GL}, in which the authors prove that
\begin{equation}\label{eq:hardyGL}
\int_{\H^d}
\frac{|z|^2}{\Kor^2} \cdot \frac{|u|^2}{\Kor^2} 
\leq
\left(\frac{2}{Q-2}\right)^2
\int_{\H^d}|\nabla_{\!\H}u|^2,
\end{equation}
for any $u\in H^1(\H^d)$.
This estimate (cf. \cite[Corollary 2.1]{GL}) represents the counterpart in the Heisenberg group of the well-known Hardy inequality \eqref{eq:hardy}, which is among the most important mathematical manifestations of the uncertainty principle from Quantum Mechanics. 

The statement of \eqref{eq:hardyGL} is quite intrinsically related to the sub-Riemannian structure of $\H^d$ (it involves a weight related to the fundamental solution of $\L$, which is given by $|(z,t)|_{\H}^{-Q+2}$, see Folland \cite{F}) and the constant of the inequality is known to be sharp. Generalizations of estimate \eqref{eq:hardyGL} are desirable, in particular in the sense of the involved weights. For our purposes, we will need other than \eqref{eq:hardyGL}, also its weighted version \eqref{eq:DamGLw} below. 
However, it is nice to show the following much more general result proved by D'Ambrosio \cite{Da2} (see also \cite{Da1, Da3} and the references therein for a comprehensive picture).
\begin{theorem}[{\cite[Theorem 3.1]{Da2}}]\label{lem:dambrosio}
	Let $d\geq1$, $\Omega\subset\H^d$ be an open set, and $h\in\mathcal{C}^1(\Omega;\R^{2d+1})$ be a vector field such that $\operatorname{div}_{\H}\, h>0$, where
	$$
	\operatorname{div}_{\H}\,h:=\operatorname{div}\,(\sigma^T\sigma h),
	\qquad
	\sigma:=
		\begin{pmatrix}
			\I_d & \0_d & 2y
			\\
			\0_d & \I_d & -2x
		\end{pmatrix}
	\in \mathcal M_{2d \times (2d+1)}(\R).
	$$
	Then the inequality
	\begin{equation}\label{eq:hardyHgeneral}
		\int_{\Omega} |u|^p \, |\operatorname{div}_{\H}\, h| \,dzdt
		\le
		p^p
		\int_{\Omega} |\sigma h|^p \, |\operatorname{div}_{\H}\, h|^{-(p-1)} |\grad u|^p \,dzdt
	\end{equation}
	holds for any $p>1$, and $u\in\mathcal{C}^1_0(\Omega)$.
\end{theorem}
In particular, setting $$h(z,t) = {|z|^{pb}}{\Kor^{-p(a+b)}} \left( z , \frac{t}{2 |z|^2} \right)$$ (to be precise, one needs to employ also a standard approximation argument to cut the singularity of $h$  when $z=0$, see \cite{Da2}), and since $|\grad\Kor| = \frac{|z|}{\Kor}$ (see \eqref{eq:kder1} and \eqref{eq:kder2} below), it follows the next
\begin{theorem}[{\cite[Theorem 3.2]{Da2}}]\label{lem:dambrosio2}
	Let $d\geq1$, $\Omega\subset\H^d$ be an open set, $p>1$, and $a,b \in \R$ such that $\frac{Q}{p} = \frac{2d+2}{p} > \max\{a,1+\frac2p-b\}$.
	Then the inequality
	\begin{equation}\label{eq:hardyH}
		\left( \frac{Q}{p}-a \right)^p
		\bigintsss_{\Omega} 
		\left[ \frac{|\grad\Kor|^b}{\Kor^a} |u| \right]^p
		dzdt
		\le
		\bigintsss_{\Omega} 
		\left[ \frac{|\grad\Kor|^{b-1}}{\Kor^{a-1}} |\grad u| \right]^p
		dzdt
	\end{equation}
	holds for $u$ in the closure of the space $\mathcal{C}^\infty_0(\Omega)$ w.r.t. the norm $\norm{ |\grad\Kor|^{b-1} \Kor^{1-a} \grad u}_{L^p(\Omega)}$.
\end{theorem}
%

Choosing in particular $\Omega=\H^d$, $d\ge1$, $p=2$, $a=b=1$, we recover the Garofalo--Lanconelli inequality \eqref{eq:hardyGL},
whereas choosing $\Omega=\H^d$, $d\ge1$, $p=2$, $a=\frac12$ and $b=1$, we obtain its weighted version
\begin{equation}\label{eq:DamGLw}
	\int_{\H^d}
	\frac{|z|^2}{\Kor^2}
	\frac{|u|^2}{\Kor}
	dzdt
	\le
	\left( \frac{2}{Q-1} \right)^2 
	\int_{\H^d}
	\Kor
	|\grad u|^2
	dzdt
	,
\end{equation}
which should be compared with its Euclidean counterpart
\begin{equation*}
	\int_{\R^d} \frac{|u|^2}{|x|} dx 
	\le
	\left( \frac{2}{d-1} \right)^2
	\int_{\R^d} |x| |\nabla u|^2 dx
	.
\end{equation*}

Even if they will not be employed in the current work, it is quite interesting to note that the result by D'Ambrosio can be used also to generate Hardy's inequalities with weights independent on the vertical direction $t$. Indeed, plugging in \eqref{eq:hardyHgeneral} the choices
$$
\Omega=\H^d,
\qquad
p=2,
\qquad
h(z,t)=\left( \frac{z}{|z|^2}, 0 \right)%
\Longrightarrow
\operatorname{div}_{\H}\,h=\frac{2d-2}{|z|^2},
\quad
|\sigma h|=\frac{1}{|z|},
$$
then \eqref{eq:hardyHgeneral} gives the horizontal Hardy inequality
\begin{equation}\label{eq:hardyDa}
\int_{\H^d}\frac{|f|^2}{|z|^2}
\leq
\left( \frac{2}{2d-2} \right)^2
\int_{\H^d}|\nabla_{\!\H}f|^2,
\end{equation}
for any $f\in H^1(\H^d)$ and $d\geq2$. 
Analogously, plugging in \eqref{eq:hardyHgeneral} the choices
$$
\Omega=\H^d,
\qquad
p=2,
\qquad
h(z,t)=\left( \frac{z}{|z|}, 0 \right)
\Longrightarrow
\operatorname{div}_{\H}\,h=\frac{2d-1}{|z|},
\quad
|\sigma h|=1,
$$
we obtain the weighted horizontal Hardy inequality
\begin{equation*}
\int_{\H^d}\frac{|f|^2}{|z|}
\leq
\left(\frac{2}{2d-1}\right)^2\int_{\H^d}|z||\nabla_{\!\H}f|^2,
\end{equation*}
for any $f\in H^1(\H^d)$, and $d\geq1$. 
Again, an approximation argument should be used to make the proof rigorous.
Compared to the Hardy inequality \eqref{eq:hardyGL} by Garofalo and Lanconelli, \eqref{eq:hardyDa} has a stronger weight at the left-hand side (obviously, $|z|^{-2}\geq|z|^2|(z,t)|_{\H}^{-4}$) but a worse constant (indeed, $\frac{2}{Q-2} = \frac1d < \frac1{d-1}$). 

Notice that in the three cases, namely Euclidean Hardy's inequality \eqref{eq:hardy}, Garofalo--Lanconelli inequality \eqref{eq:hardyGL} and its horizontal version \eqref{eq:hardyDa} by D'Ambrosio, the best constant of the inequality has the usual shape $\big(\frac{2}{D-2}\big)^2$, where the dimensional role $D$ is played respectively by the Euclidean dimension $d$ in \eqref{eq:hardy}, by the \emph{homogenous dimension} $Q$ in \eqref{eq:hardyGL}, and by the \emph{horizontal dimension} $2d$ in the case of \eqref{eq:hardyDa}. Similarly in their weighted version.

For completeness, we cite the recent developments in \cite{CFKP} about magnetic Hardy inequalities in the context of the Heisenberg group. 

With \eqref{eq:hardyGL} and \eqref{eq:DamGLw} in our hands, we are ready to prove our main results.

\section{Proof of Theorem\til\ref{thm:grad-est}}\label{sec:proof}

In the sequel, by the identification $\C^d \simeq \R^{2d}$, the symbol $z$ will be used for the vector $(x,y) \in \R^{2d}$. Also, for short we will always write the integrals 
\begin{equation*}
	\int f = \int_{\H^d}f(z,t) \,dz\,dt,
	\qquad
	f:\H^d\to\C
	,
\end{equation*}
and we will make use of the notation
\begin{equation}\label{eq:Korpesos}
	\kor := \Kor,
	\qquad
	\gl := \frac{|z|}{\Kor}
	.
\end{equation}
It is straightforward to compute
%
%
\begin{equation}\label{eq:kder1}
	\grad |(z,t)|_{\H}
	=
	\frac{|z|^2 z + (y,-x) t}{|(z,t)|_{\H}^{3}},
	\qquad
	\L |(z,t)|_{\H} 
	=
	-
	(2d+1) \frac{|z|^2}{|(z,t)|_{\H}^3},
\end{equation}
and hence
\begin{equation}\label{eq:kder2}
	|\grad\kor| = \gl,
	\qquad
	\L\kor = - (Q-1) \frac{\gl^2}{\kor},
\end{equation}
which should be compared to the Euclidean case $|\nabla|x||=1$, $-\Delta|x|=-(d-1)\frac{1}{|x|}$.

\subsection{Proof of estimate \eqref{eq:est1}}\label{subsec:est1}
The case $|\lambda_2|> \delta \lambda_1$ splits into the two subcases
$\l_2> \delta \l_1$ and $\l_2<- \delta \l_1$. Let us firstly consider the case $\l_2>\delta\l_1$. We may assume that $u \in H^1(\H^d)$ is a non-null solution of \eqref{eq:prob0} (otherwise, estimate \eqref{eq:est1} is trivial). Let us choose $v=u$ in \eqref{eq:solution} as test function:
subtracting $\delta$ times the real part of the resulting identity from its imaginary part,
we get
	\begin{equation}\label{eq:|l2|>l1-1}
		(\delta\l_1-\l_2) \int |u|^2 = \delta \int |\grad u|^2 + \delta \re \int \o{u} f - \im \int \o{u} f.
	\end{equation}
	By Cauchy--Schwarz and the Hardy inequality \eqref{eq:hardyGL}, we estimate
	\begin{multline*}
		\left|
			\delta \re \int \o{u} f - \im \int \o{u} f
		\right|
		\le
		(\delta+1) \int |u||f|
		\le
		(\delta+1) \norm{ \frac{\kor}{\gl} f} \norm{ \frac{\gl}{\kor} u}
		\le
		\frac{\delta+1}{d} \norm{ \frac{\kor}{\gl} f} \norm{\grad u}.
	\end{multline*}	
	Plugging the last inequality in \eqref{eq:|l2|>l1-1}, we get
	\begin{equation}\label{eq:|l2|>l1-12}
	0 > (\delta\l_1-\l_2) \norm{u}^2
	\ge 
	\delta \norm{\grad u} \left(\norm{\grad u} - \frac{1+1/\delta}{d} \norm{ \frac{\kor}{\gl} f} \right)
	\end{equation}
	from which, since $u$ is non-vanishing, estimate \eqref{eq:est1} necessarily follows.
	
	In the case $\l_2<- \delta \l_1$, after setting $v=u$ in \eqref{eq:solution} and taking the imaginary part and $\delta$ times the real part, we add the two equations obtaining 
	\begin{equation*}
		(\delta \l_1+\l_2) \int |u|^2 = \delta \int |\grad u|^2 + \delta \re \int \o{u} f + \im \int \o{u} f
	\end{equation*}
	in place of \eqref{eq:|l2|>l1-1}. Then the proof proceeds in a completely analogous way as above.

\subsection{Proof of estimate \eqref{eq:est2}}
We assume now $|\lambda_2|\leq \delta \lambda_1$, and
notice that then $\lambda_1\geq0$. 
The proof, which is in part inspired by \cite{FKV}, is based on algebraic manipulations of equation \eqref{eq:prob0} (or more precisely, of its weak formulation \eqref{eq:solution}), and suitable choices of the test functions. 
For the sake of clarity, we first give the proof assuming $u,f$ to be smooth enough and compactly support; in a final step, we will remove these regularity assumptions.

\subsubsection{The three identities}

Let $\Phi_1, \Phi_2, \Phi_{3,1}, \Phi_{3,2} \colon \H^d \to \R$ be four sufficiently smooth functions. 
Choosing $v=\Phi_1 u$ in equation \eqref{eq:solution}, and taking the real part of the resulting identity, we get
	\begin{equation}\label{eq:fond1}
		- \frac{1}{2} \int \L\Phi_1 |u|^2 
		- \int \Phi_1 |\grad u|^2
		+
		\lambda_1 \int \Phi_1 |u|^2
		=
		\re \int f \Phi_1 \o{u}
	\end{equation}
	where we used the identity 
	\begin{equation}\label{eq:gradu2}
		2 \, \re(\o{u} \grad u) = \grad |u|^2
	\end{equation}
	and the integration by parts in the first term.
Similarly, choosing $v=\Phi_2 u$ in equation \eqref{eq:solution}, and taking the imaginary part, we obtain the identity
	\begin{equation}\label{eq:fond2}
		-\im \int \grad \Phi_2 \cdot \o{u}\grad u
		+
		\lambda_2 \int \Phi_2 |u|^2 
		=
		\im \int f \Phi_2 \o{u} .
	\end{equation}
Finally, let us consider the skew-symmetric multiplier
\begin{equation}\label{def:skewmult}
		v = [\L,\Phi_{3,1}] u - [T^2, \Phi_{3,2}] u
		= u \L \Phi_{3,1} - 2\grad \Phi_{3,1} \cdot \grad u - uT^2\Phi_{3,2} - 2 T\Phi_{3,2} \, Tu
\end{equation}
where $T^2 = TT$. We are inspired by the Euclidean case, where one chooses the multiplier in the form $v=[-\Delta,\Phi_3]u$, with $\Phi_3(x)=\frac{|x|^2}{2}$. On one hand we have that
\begin{equation}\label{eq:3.1}
		\begin{split}
			-\re \int \grad([\L, \Phi_{3,1}] \o{u}) \cdot \grad u
			=&\,
			-\frac{1}{2} \int \L^2 \Phi_{3,1} |u|^2
			+
			2 \int \grad \o{u} \cdot \grad^2 \Phi_{3,1} \grad u
			\\
			&
			+ 8 \re \int \grad\Phi_{3,1} \cdot \begin{pmatrix} \0_d & \I_d \\ -\I_d & \0_d \end{pmatrix} \grad\o{u} Tu 
		\end{split}
\end{equation}
where $\L^2=\L\L$ is the bi-sublaplacian, $\0_d$ and $\I_d$ are the $d$-dimensional zero and identity matrix respectively, $\grad^2$ is the Hessian matrix 
\begin{equation*}
	\grad^2 =
	\begin{pmatrix}
		X_i X_j & X_i Y_j
		\\
		Y_i X_j & Y_i Y_j
	\end{pmatrix}_{i,j \in \{1,\dots,d\}}
\end{equation*}
and where we used the integration by parts combined with identity \eqref{eq:gradu2} and its analogue $$2\re(\grad^2 u \grad \o{u}) = \grad |\grad u|^2.$$ 
Similarly we also have 
\begin{equation}\label{eq:3.2}
		\re \int \grad([T^2, \Phi_{3,2}] \o{u}) \cdot \grad u
		=
		\frac{1}{2} \int \L T^2 \Phi_{3,2} |u|^2
		+
		2
		\re \int \grad T\Phi_{3,2} \cdot \grad\o{u} Tu 
		.
\end{equation}
On the other hand
\begin{equation}
	\label{eq:3.3}
	\re \int \lambda ([\L, \Phi_{3,1}] \o{u}) u
	=
	-
	2\lambda_2 \im \int \grad\Phi_{3,1} \cdot \o{u} \grad u
\end{equation}
using again \eqref{eq:gradu2}, the fact that $\re(ab)= \re a\re b - \im a \im b$ for any $a,b \in \C$ and that $\Phi_{3,1}$ is real-valued, and similarly
\begin{equation}
	\label{eq:3.4}
	\re \int \lambda ([T^2, \Phi_{3,2}] \o{u}) u
	=
	2\lambda_2 \im \int T\Phi_{3,2} \cdot \o{u} T u
	.
\end{equation}

	\begin{remark}
	Compared to the Euclidean case, from a technical point of view, the computations here are complicated  by the fact that the Hessian $\grad^2$ is not symmetric, due to the commutation relations \eqref{eq:commutators}, from which it originates the last term on the right-hand side of \eqref{eq:3.1}.
	\end{remark}
	
Choosing the multiplier \eqref{def:skewmult} in \eqref{eq:solution}, considering the real part of the resulting identity multiplied by a factor $1/2$, and taking into account \eqref{eq:3.1}, \eqref{eq:3.2}, \eqref{eq:3.3} and \eqref{eq:3.4} we finally obtain the identity
%
%
	\begin{equation}\label{eq:fond3}
		\begin{split}
			&
			\re \int \grad \o{u} \cdot \grad^2 \Phi_{3,1} \grad u
			-
			\frac{1}{4} \int \left[ \L^2 \Phi_{3,1} - \L T^2 \Phi_{3,2} \right] |u|^2
			\\
			&
			+ \re \int 
			\left[
			4 
			\grad\Phi_{3,1} \cdot
			\begin{pmatrix} \0_d & \I_d \\ -\I_d & \0_d \end{pmatrix}  
			+
			\grad T\Phi_{3,2}
			\right]
			\cdot
			\grad\o{u} Tu 
			\\
			&
			-
			\lambda_2 \im \int \grad\Phi_{3,1} \cdot \o{u} \grad u
			-
			\lambda_2 \im \int T\Phi_{3,2} \cdot \o{u} T u
			\\
			=&\,
			\frac{1}{2} \re \int f \left[ \L \Phi_{3,1}-T^2\Phi_{3,2} \right] \o{u}  
			- \re \int f \grad\Phi_{3,1} \cdot \grad\o{u}
			- \re \int f T\Phi_{3,2} \cdot T\o{u}
			.
		\end{split}
	\end{equation}
Roughly speaking, the idea of the multipliers method we are here employing is to recognize in the equation \eqref{eq:fond3} part of the square $|\grad (e^{-i\psi} u)|^2$ for some scalar function $\psi$, and to complete it adding the symmetric identities with suitable choices of $\Phi_1$ and $\Phi_2$. At this aim, first of all we would like to have $\grad^2 \Phi_{3,1} = \I_{2d}$ in the first line of \eqref{eq:fond3}, and also to kill the bad term in the second line; therefore, we choose
\begin{equation*}
	\Phi_{3,1}(z,t) = \frac{|z|^2}{2},
	\qquad
	\Phi_{3,2}(z,t) = t^2,
\end{equation*}
so that our third inequality reads as
\begin{equation}\label{eq:fond3.0}
	\int |\grad u|^2 - \lambda_2 \im \int \o{u} \left[ z \cdot \nabla_z u + 2t Tu \right]
	= 
	- (d+1) \re \int \o{f} u - \re \int \o{f} \left[ z \cdot \nabla_z u + 2t Tu \right] .
\end{equation}

\begin{remark}
	Other than mere technical reason, there is a much more relevant and intrinsic argument supporting the choice of the multiplier
	\begin{equation}\label{eq:mult3}
		v = \left[\L,\frac{|z|^2}{2} \right] u - [T^2,t^2]u 
	\end{equation}
	as the adequate one,
	relying on the strong connection between the multipliers method and the Mourre theory, as firstly noticed in \cite{CK} and reasoned in \cite{CFS}. By analogy, let us firstly consider the Euclidean case: here, the antisymmetric multiplier actually is, up to a constant factor, the generator of the (classic Euclidean) dilations, namely
	\begin{equation*}
		v = \left[ -\Delta, \frac{|x|^2}{2} \right]u = -i \left( x \cdot \nabla_x + \frac{d}{2} \right) u =: -2i Au 
	\end{equation*}
	where the dilation acts like $\widetilde{\delta}_\lambda f(x)= e^{\frac{d}{2} \lambda} f(e^\lambda x)$ and $\widetilde{\delta}_\lambda = e^{i\lambda A}$.
	
	When we turn to the Heisenberg case, we need to substitute the Euclidean dilation with the Heisenberg dilations, namely
	\begin{equation*}
		\delta_\lambda f(z,t)= e^{\frac{Q}{2} \lambda} f(e^\lambda z, e^{2\lambda} t).
	\end{equation*} 
	The generator of dilations ($\delta_\lambda = e^{i\lambda A_{\H}}$) is now given by
	\begin{equation*}
		A_{\H} := -i \left( z\cdot\nabla_{\!z} + 2t T + \frac{Q}{2} \right) ,
	\end{equation*}
	and, in the light of this reasoning, is it not surprising now to notice that the multiplier \eqref{eq:mult3} is actually equal to $v=-2i A_{\H} u$. 
\end{remark}

With the identity \eqref{eq:fond3.0} in our hands, to continue our argument we would like to write
\begin{equation}\label{eq:cond1}
	z \cdot \nabla_z u + 2t Tu 
	=
	C(z,t) \grad\psi \cdot \grad u
\end{equation}
for some non-null functions $C$ and $\psi$. However, this is not possible without any assumption on $u$. Indeed, it is readily checked that
\begin{equation*}
	\grad \psi \cdot \grad u = \nabla_z \psi \cdot \nabla_z u + (y,-x) T\psi \nabla_z u + (y,-x) \nabla_z \psi Tu + |z|^2 T\psi Tu
	,
\end{equation*}
and imposing \eqref{eq:cond1} we get that $\psi$ should satisfy
\begin{equation*}
	\nabla_z + (y,-x) T\psi = Cz,
	\qquad
	(y,-x) \cdot \nabla_z \psi + |z|^2 T\psi = C 2t.
\end{equation*}
Multiplying the first one times $(y,-x)$ and using the second one, one gets $2t C =0$: a contradiction.
We will assume $u$ to be radial respect to the Koranyi norm, motivated also by the fact that 
\begin{equation*}
	z \cdot \nabla_z + 2t T = \mathbb{E} = \kor \, \partial_{\kor}
	,
\end{equation*}
where $\mathbb{E}$ is the Euler vector field (recall \eqref{eq:EulerVF}) and $\partial_{\kor}$ is the radial derivative with respect to the Koranyi norm in \eqref{eq:koranyi}.  Here is also the key to choose this kind of radial functios rather than functions whose radiality is with respect to the variable $|z|$.
However, before doing this, it will be convenient to introduce a suitable change of variables, and rewrite our identities accordingly.


\subsubsection{The change of variables}
\label{sub:change}

Let us consider the change of variables $\phi\colon \H^d \setminus \mathcal{Z} \to \H^d \setminus \mathcal{Z}$:
\begin{align*}
	(z,t) &= \phi(w,s) = \left( \frac{w}{\jap{s}^{1/2}} , \frac{s}{\jap{s}}|w|^2 \right) ,
	\\
	(w,s) &= \phi^{-1}(z,t) = \left( \frac{|(z,t)|_{\H}}{|z|} z , \frac{t}{|z|^2} \right) ,
\end{align*}
where $\mathcal{Z} := \{ (z,t) \in \H^d \colon z=0 \}$ is the center of the Heisenberg group, and $\jap{s}=\sqrt{1+s^2}$ is the Japanese brackets notation. This change of variables satisfy the following identities:
\begin{equation}\label{hats}
	|w|=|(z,t)|_{\H}=\kor,
	\qquad
	\frac{1}{\jap{s}}=\frac{|z|^2}{|(z,t)|_{\H}^2}=\gl^2,
	\qquad
	\hat{w} := \frac{w}{|w|}=\frac{z}{|z|} =: \hat{z},
\end{equation}
where $\kor$ and $\gl$ are defined in \eqref{eq:Korpesos}.
This means that the unit Koranyi sphere punctured at the north and south poles, viz. $\{ (z,t) \in \H^d \setminus \mathcal{Z} \colon \kor=1 \}$, is mapped by $\phi^{-1}$ in the unit cylinder with respect to the new coordinates, viz. $\{ (w,s) \in \H^d \setminus \mathcal{Z} \colon |w|=1 \}$; as a consequence, a radial function $f(z,t)=f^0(\Kor)$ is a cylindrical function in the new coordinates, in the sense that $(f\circ\phi)(w,s)=f^0(|w|)$.
Moreover, the Garofalo--Lanconelli weight appearing in \eqref{eq:hardyGL} now depends only on the new vertical variable $s$.

The Jacobian matrix of this transformation is given by
\begin{equation*}
	\J\phi
	=
	\begin{pmatrix}
		\dfrac{1}{\jap{s}^{1/2}} \I_{2d} & -\dfrac12 \dfrac{s}{\jap{s}^{5/2}} w
		\\
		2\dfrac{s}{\jap{s}}w & \dfrac{|w|^2}{\jap{s}^3}
	\end{pmatrix}
	=
	\begin{pmatrix}
		\dfrac{|z|}{\Kor} \I_{2d} & -\dfrac12 \dfrac{|z|^3}{\Kor^4} t\hat{z}
		\\
		2\dfrac{t}{\Kor} \hat{z} & \dfrac{|z|^6}{\Kor^4}
	\end{pmatrix}
\end{equation*}
and hence
\begin{equation*}
	|\J_\phi| := |\det \J\phi| = \frac{|w|^2}{\jap{s}^{d+1}},
	\qquad
	\int_{\H^d} f(z,t) dzdt = \int_{\R_*^{2d}\times \R} |\J_{\phi}| f(\phi(w,s)) dwds.
\end{equation*}
%
%
On the other side, the Jacobian matrix of the inverse transformation is given by
\begin{align*}
	\J\phi^{-1}
	=
	\begin{pmatrix}
		\dfrac{\Kor}{|z|} \left[ \I_{2d} - \dfrac{t^2}{\Kor^4} \hat{z}\otimes\hat{z} \right]
		&
		\dfrac12 \dfrac{t\hat{z}}{\Kor^3}
		\\
		-2 \dfrac{t\hat{z}}{|z|^3}
		&
		\dfrac{1}{|z|^2}
	\end{pmatrix}
	=
	\begin{pmatrix}
		\jap{s}^{\frac12} \left[ \I_{2d} - \dfrac{s^2}{\jap{s}^2} \hat{w}\otimes\hat{w} \right] & \dfrac12 \dfrac{s}{\jap{s}} \dfrac{\hat{w}}{|w|} 
		\\
		-2 s \jap{s}^{\frac12} \dfrac{\hat{w}}{|w|} & \dfrac{\jap{s}}{|w|^2}
	\end{pmatrix}.
\end{align*}

Now, let us set $\u := u \circ \phi$, so that
$
	\u(w,s) = u(\phi(w,s)) = u(z,t)
$
. 
By the chain rule
%
\begin{equation*}
	(\nabla_z, \partial_t) u(z,t)
	=
	(\nabla_w, \partial_s) \u(w,s) \cdot \J\phi^{-1}
\end{equation*}
and hence
\begin{align*}
	\nabla_z u(z,t) &= \jap{s}^{1/2} \left[ \I_{2d} - \dfrac{s^2}{\jap{s}^2} \hat{w}\otimes\hat{w} \right] \cdot \nabla_w \u(w,s) 
	- 2s\jap{s}^{1/2} \frac{\hat{w}}{|w|} \partial_s \u(w,s)
	\\
	\partial_t u(z,t) &= \frac12 \frac{s}{\jap{s}} \frac{\hat{w}}{|w|} \cdot \nabla_w \u(w,s) 
	+  \dfrac{\jap{s}}{|w|^2} \partial_s \u(w,s) .
\end{align*}
Denoting for short 
\begin{equation}\label{eq:stars}
	w^\perp := J w,
	\qquad
	z^\perp := J z = (y,-x) ,
	\qquad
	J := \begin{pmatrix} \0_d & \I_d \\ -\I_d & \0_d \end{pmatrix},
\end{equation}
we have
\begin{equation*}
	\begin{split}
		\grad u(z,t) 
		=\,&\, 
		\nabla_z u(z,t) + 2 z^\perp Tu(z,t)
		\\
		=\,&\,
		\jap{s}^{1/2} \nabla_w \u(w,s)
		-
		\jap{s}^{1/2}
		\left[ \frac{s}{\jap{s}} \hat{w} - \frac{1}{\jap{s}} \hat{w}^\perp \right]
		\left[ \frac{s}{\jap{s}} \hat{w} \cdot \nabla_w \u(w,s) + 2 \frac{\jap{s}}{|w|} \partial_s  \u(w,s) \right]
		\\
		=:&\,
		\wgrad \u(w,s)
		.
	\end{split}
\end{equation*}

\begin{remark}\label{rmk:spherical}
	Another meaningful change of variables, sometimes employed in the literature (see e.g. \cite{Greiner,GK,Da0}), is the one involving some modified spherical coordinates, namely
	\begin{equation*}
		(z,t) = \widetilde{\phi}(\rho,\theta,\vartheta) = (\rho \sqrt{\sin\theta} \, \omega(\vartheta), \rho^2 \cos\theta)
	\end{equation*}
	where the unit vector $\omega(\vartheta)\in\R^{2d}$ is given by
	\begin{equation*}
		\omega(\vartheta) = ( \cos\vartheta_1 , \sin\vartheta_1 \cos\vartheta_2 , \dots, \sin\vartheta_1 \sin\vartheta_2 \cdots \cos\vartheta_{2d-1} , \sin\vartheta_1 \sin\vartheta_2\cdots \sin\vartheta_{2d-1})
	\end{equation*}
	and $\rho\ge0$, $\theta \in [0,\pi]$, $\vartheta=(\vartheta_1,\dots,\vartheta_{2d-2},\vartheta_{2d-1}) \in [0,\pi] \times \cdots \times [0,\pi] \times [0,2\pi[$.
	
	Actually, we are essentially using this change of variables: indeed it is easy to switch from one to the other observing that 
	\begin{equation*}
		\rho = |w| , 
		\quad
		\omega(\vartheta) = \hat{w} ,
		\quad
		\theta= \arccos\left(\frac{s}{\jap{s}}\right) ,
		\quad
		\frac{1}{\jap{s}} = \sin\theta,
		\quad
		\frac{s}{\jap{s}} = \cos\theta
		.
	\end{equation*}
	Moreover, the unitary vector $\jap{s}^{-1} \left[ s \hat{w} - \hat{w}^\perp \right]$ appearing in $\wgrad$ can be rewritten as the action of a rotation matrix on $\hat{w}$, viz.
	\begin{equation*}
		\frac{s}{\jap{s}} \hat{w} - \frac{1}{\jap{s}} \hat{w}^\perp = R_\theta \hat{w}
		,
		\qquad
		R_\theta 
		:=
		\begin{pmatrix}
			\cos\theta \, \I_d & -\sin\theta \, \I_d \\ \sin\theta \, \I_d & \cos\theta \, \I_d
		\end{pmatrix}
		.
	\end{equation*}

	Using the spherical change of variables $\widetilde{\phi}$ it is also immediate to see that
	\begin{equation*}
		\rho \, \partial_\kor  
		= \rho \frac{\partial}{\partial\kor}(z,t) \cdot (\nabla_z,\partial_t) 
		= \rho (\sqrt{\sin\theta} \, \omega(\vartheta), 2\rho \cos\theta) \cdot (\nabla_z,\partial_t) 
		= (z,2t) \cdot (\nabla_z,\partial_t)
		= \mathbb{E} 
	\end{equation*}
	justifying the connection between the Euler vector field and the radial derivative with respect to the Koranyi gauge.
\end{remark}

Exploiting the change of variables $(z,t)=\phi(w,t)$, the three identities \eqref{eq:fond1}, \eqref{eq:fond2} and \eqref{eq:fond3.0} becomes respectively
	\begin{gather}
	\label{eq:fondws1}
		- \frac{1}{2} \int|\J_\phi| \widetilde{\L} \Psi_1 |\u|^2 
		- \int|\J_\phi| \Psi_1 |\wgrad \u|^2
		+
		\lambda_1 \int|\J_\phi| \Psi_1 |\u|^2
		=
		\re \int|\J_\phi| \f \Psi_1 \o{\u}
		,
	\\
	\label{eq:fondws2}
		-\im \int|\J_\phi| \wgrad \Psi_2 \cdot \o{\u}\wgrad \u
		+
		\lambda_2 \int|\J_\phi| \Psi_2 |\u|^2 
		=
		\im \int|\J_\phi| \f \Psi_2 \o{\u} 
		,
	\\
	\label{eq:fondws3}
		\int |\J_\phi| |\wgrad \u|^2 - \lambda_2 \im \int |\J_\phi| w \cdot \o{\u} \nabla_w \u  
		= 
		- (d+1) \re \int |\J_\phi| \o{\f} \u - \re \int |\J_\phi| \o{\f} \, w \cdot \nabla_w \u  
		,
	\end{gather}
where we set $\f = f \circ \phi$, $\Psi_1 = \Phi_1 \circ \phi $, $\Psi_2 = \Phi_2 \circ \phi$ and $\widetilde{\L}$ is $\L$ in the new variables ($\widetilde{\L}\Psi_1:=\L\Phi_1$). Notice in particular that
\begin{equation}\label{eq:radialder}
	w \cdot \nabla_w \u(w,s) = (z \cdot \nabla_z + 2t \partial_t) u(z,t) 
\end{equation}
and that $\hat{w} \cdot \nabla_w$ is the radial derivative with respect to the Koranyi norm in these coordinates.

\subsubsection{Radial assumptions}
It is time to finally use our assumptions of radiality for the solutions to proceed in the argument. Firstly notice that, for any radial function $g(w,s) = g^0(|w|)$ it holds
\begin{gather*}
	\wgrad g = \frac{1}{\jap{s}^{1/2}} \left[ \frac{1}{\jap{s}} \hat{w} + \frac{s}{\jap{s}} \hat{w}^\perp \right] g^0_r
	,
	\qquad
	|\wgrad g| = \frac{|g^0_r|}{\jap{s}^{1/2}} 
	,
	\qquad
	\widetilde{\L} g = -\frac{1}{\jap{s}} \left( g^0_{rr} + \frac{2d+1}{|w|} g^0_r \right)
	.
\end{gather*}
Assuming $\u(w,s)=\u^0(r)$, $\Psi_1(w,s)=\Psi_1^0(r)$, and $\Psi_2(w,s)=\Psi_2^0(r)$ with $r=|w|$, the identities \eqref{eq:fondws1}, \eqref{eq:fondws2} and \eqref{eq:fondws3} become respectively
	\begin{gather*}
		\frac{1}{2} \int \frac{|\J_\phi|}{\jap{s}} \left( (\Psi_1^0)_{rr} + \frac{2d+1}{|w|} (\Psi_1^0)_r \right) |\u^0|^2 
		- \int \frac{|\J_\phi|}{\jap{s}} \Psi_1^0 |\u^0_r|^2
		+
		\lambda_1 \int|\J_\phi| \Psi_1^0 |\u^0|^2
		=
		\re \int|\J_\phi| \f \Psi_1^0 \o{\u^0} ,
	\\
		-\im \int \frac{|\J_\phi|}{\jap{s}} (\Psi_2^0)_r \cdot \o{\u^0} \u^0_r
		+
		\lambda_2 \int|\J_\phi| \Psi_2^0 |\u^0|^2 
		=
		\im \int|\J_\phi| \f \Psi_2^0 \o{\u^0} ,
	\\
		\int \frac{|\J_\phi|}{\jap{s}} |\u^0_r|^2 
		- \lambda_2 \im \int |\J_\phi| |w| \cdot \o{\u^0} \u^0_r  
		= 
		- (d+1) \re \int |\J_\phi| \o{\f} \u^0 - \re \int |\J_\phi| \o{\f} \, |w| \u^0_r
		.
	\end{gather*}

Let us write the following elementary computation as a lemma.
\begin{lemma}\label{lem:B}
	Let $\alpha>0$. Then
	\begin{equation*}
		\int_{-\infty}^{+\infty} \frac{ds}{\jap{s}^{\alpha+1}} 
		=
		\sqrt{\pi} \frac{ \Gamma\left( \frac{\alpha}{2} \right)}{\Gamma\left( \frac{\alpha+1}{2} \right)},
	\end{equation*}
	In particular, 
	by Fubini--Tonelli Theorem, it follows
	\begin{equation*}
		\int_{\R^d \times \R}
		\frac{|f(w)|}{\jap{s}^{\gamma+d+1}}
		dwds
		=
		B_{\gamma+d,d} 
		\int_{\R^d \times \R}
		\frac{|f(w)|}{\jap{s}^{d+1}}
		dwds
		,
		\qquad
		B_{\alpha,\beta} 
		:= 
		\frac{\Gamma(\frac{\alpha}{2}) \Gamma(\frac{\beta+1}{2})}{\Gamma(\frac{\beta}{2}) \Gamma(\frac{\alpha+1}{2})} 
		=
		\frac{1}{B_{\beta,\alpha}}
		,
	\end{equation*}
	for any $\gamma+d>0$ and $f$ depending only on $w$.
\end{lemma}
To calculate the integral in Lemma \ref{lem:B}, we just make change of variable $r=(1+s^2)^{-1}$, which implies $\int_\R\jap{s}^{-\alpha-1}ds=B(\frac\alpha 2,\frac12)=\Gamma(\frac \alpha 2)\Gamma(\frac 12)/\Gamma(\frac{\alpha+1}{2})$, where $B(a,b)=\int_0^1 r^{a-1}(1-r)^{b-1}dr$ is the beta function.

Recalling that $|\J_\phi| = \frac{|w|^2}{\jap{s}^{d+1}}$, using the lemma with $\gamma=1$, and defining for short 
\begin{equation}\label{eq:Gd}
	G_d := B_{d+1,d} = \frac2d \left[ \frac{\Gamma\left(\frac{d+1}{2}\right)}{\Gamma\left(\frac{d}{2}\right)} \right]^2 ,
\end{equation}
the three identities are equivalent to
	\begin{align*}
			\frac{1}{2} \int |\J_\phi| \left( (\Psi_1^0)_{rr} + \frac{2d+1}{|w|} (\Psi_1^0)_r \right)  |\u^0|^2 
			&- \int |\J_\phi| \Psi_1^0 |\u^0_r|^2
			+
			\frac{\lambda_1}{G_d} \int|\J_\phi| \Psi_1^0 |\u^0|^2
			\\
			&=
			\frac{1}{G_d}
			\re \int|\J_\phi| \o{\f} \Psi_1^0 \u^0 ,
	\\
		-\im \int |\J_\phi| (\Psi_2^0)_r \cdot \o{\u^0} \u^0_r
		+
		\frac{\lambda_2}{G_d} \int|\J_\phi| \Psi_2^0 |\u^0|^2 
		&=
		- \frac{1}{G_d}
		\im \int|\J_\phi| \o{\f} \Psi_2^0 \u^0 ,
	\\
		2\int |\J_\phi| |\u^0_r|^2 
		- 2\frac{\lambda_2}{G_d} \im \int |\J_\phi| |w| \cdot \o{\u^0} \u^0_r  
		&= 
		- 2 \frac{d+1}{G_d} \re \int |\J_\phi| \o{\f} \u^0 
		- 2 \frac{1}{G_d} \re \int |\J_\phi| \o{\f} \, |w| \u^0_r
		.
	\end{align*}

Now, recalling \eqref{eq:upm}, we have that
\begin{equation*}
	\u^-(w,s) := (u^- \circ \phi)(w,s) = e^{-i\widetilde{\psi}(w)} \u(w,s)
	,
	\qquad
	\widetilde{\psi}(w) = \sgn(\lambda_2) \sqrt{\tfrac{\lambda_1}{G_d}} \, |w|
	.
\end{equation*}
Summing up all the three identities and choosing
\begin{equation*}
	\Psi^0_1(|w|) = 1 - \frac{|\lambda_2/G_d|}{\sqrt{\lambda_1/G_d}} |w| ,
	\qquad
	\Psi^0_2(|w|) = 2\widetilde{\psi}(w) = 2\sgn(\lambda_2) \sqrt{\lambda_1/G_d} \, |w| ,
\end{equation*}
we obtain
\begin{align*}
	&\int |\J_\phi| |(e^{-i\widetilde{\psi}}\u^0)_r|^2 \left[ 1 + \frac{|\lambda_2/G_d|}{\sqrt{\lambda_1/G_d}}|w| \right] 
	-
	\frac{2d+1}{2}
	\frac{|\lambda_2/G_d|}{\sqrt{\lambda_1/G_d}}
	\int |\J_\phi| \frac{|\u^0|^2}{|w|}
	\\
	=&
	- \frac{1}{G_d} \im \int |\J_\phi| 2 \sgn(\lambda_2) \sqrt{\lambda_1/G_d} |w| \o{\f} \u^0
	- \frac{1}{G_d} \re \int |\J_\phi| \o{\f} \left[(2d+1) \u^0 + \frac{|\lambda_2/G_d|}{\sqrt{\lambda_1/G_d}}|w| \u^0 + 2 |w| \u^0_r\right]
	\\
	=&
	- \frac{1}{G_d} \re \int |\J_\phi| \o{\f} \left[(2d+1) \u^0 + \frac{|\lambda_2/G_d|}{\sqrt{\lambda_1/G_d}}|w| \u^0 + 2 |w| e^{i\widetilde{\psi}} (e^{-i\widetilde{\psi}} \u^0)_r \right]
	,
\end{align*}
where we used that
\begin{align*}
	(e^{-i\widetilde{\psi}} \u^0)_r &= e^{-i\widetilde{\psi}} \left[ \u^0_r - i \sgn(\lambda_2) \sqrt{\tfrac{\lambda_1}{G_d}} \u^0 \right] ,
	\\
	|(e^{-i\widetilde{\psi}} \u^0)_r|^2 &= |\u^0_r|^2 + \tfrac{\lambda_1}{G_d} |\u^0|^2 - \im\left( 2 \sgn(\lambda_2) \sqrt{\tfrac{\lambda_1}{G_d}} \, \o{\u^0} \u^0_r \right).
\end{align*}
Multiplying by $G_d$ and using again Lemma\til\ref{lem:B} with $\gamma=1$, we arrive at
\begin{multline*}
	\int |\J_\phi| \frac{|(e^{-i\widetilde{\psi}}\u^0)_r|^2}{\jap{s}} \left[ 1 + \frac{|\lambda_2/G_d|}{\sqrt{\lambda_1/G_d}}|w| \right] 
	-
	\frac{2d+1}{2}
	\frac{|\lambda_2/G_d|}{\sqrt{\lambda_1/G_d}}
	\int |\J_\phi| \frac{1}{\jap{s}} \frac{|\u^0|^2}{|w|}
	\\
	=
	- \re \int |\J_\phi| \o{\f} \left[(2d+1) \u^0 + \frac{|\lambda_2/G_d|}{\sqrt{\lambda_1/G_d}}|w| \u^0 + 2 |w| e^{i\widetilde{\psi}} (e^{-i\widetilde{\psi}} \u^0)_r \right]
	.
\end{multline*}
Using the formulas
\begin{equation*}
	|\wgrad g|^2 = \frac{1}{\jap{s}} |g^0_r|^2 ,
	\qquad
	\jap{s}^{1/2} \left[ \frac{1}{\jap{s}} w + \frac{s}{\jap{s}} w^\perp \right] \cdot \wgrad g = |w| g^0_r ,
\end{equation*}
valid for any $g(w,s) = g^0(|w|)$, we finally reach the key identity
\begin{multline*}
		 \int |\J_\phi| |\wgrad \u^-|^2 \left[ 1 + \frac{|\lambda_2/G_d|}{\sqrt{\lambda_1/G_d}}|w| \right] 
		-
		\frac{Q-1}{2}
		\frac{|\lambda_2/G_d|}{\sqrt{\lambda_1/G_d}}
		\int |\J_\phi| \frac{1}{\jap{s}} \frac{|\u^-|^2}{|w|}
		\\
		=	
		-
		\re \int |\J_\phi| \jap{s}^{1/2} |w| \o{e^{-i\widetilde{\psi}} \f} 
		\left[ \frac{Q-1}{\jap{s}^{1/2}} \frac{\u^-}{|w|}
		+ 
		\frac{|\lambda_2/G_d|}{\sqrt{\lambda_1/G_d}} \frac{1}{\jap{s}^{1/2}} \u^- 
		+ 
	2 \left[ \frac{1}{\jap{s}} \hat{w} + \frac{s}{\jap{s}} \hat{w}^\perp \right] \cdot \wgrad \u^- \right].
\end{multline*}
Changing back the variables, and recalling the definitions of $\kor,\gl$ in \eqref{eq:Korpesos}, $\hat{w}, \hat{z}$ in \eqref{hats}  and $w^\perp,z^\perp$ in \eqref{eq:stars}, this is equivalent to 
\begin{equation}\label{eq:key2}
	\begin{split}
		& 
		\int |\grad u^-|^2 \left[ 1 + \frac{|\lambda_2/G_d|}{\sqrt{\lambda_1/G_d}} \kor \right] 
		-
		\frac{Q-1}{2}
		\frac{|\lambda_2/G_d|}{\sqrt{\lambda_1/G_d}}
		\int \gl^2 \frac{|u^-|^2}{\kor}
		\\
		&=
		-
		\re \int \frac{\kor}{\gl} \o{e^{-i\psi} f} 
		\left[ (Q-1) \gl \frac{u^-}{\kor}
		+ 
		\frac{|\lambda_2/G_d|}{\sqrt{\lambda_1/G_d}} \gl u^-
		+ 
		2 \left[ \frac{|z|^2}{\kor^2} \hat{z} + \frac{t}{\kor^2} \hat{z}^\perp \right] \cdot \grad u^- \right]
		.
	\end{split}
\end{equation}

\subsubsection{Estimates}

We are ready to start our estimates, employing the Garofalo--Lanconelli estimate \eqref{eq:hardyGL} and its weighted version \eqref{eq:DamGLw}.
For the left-hand side of \eqref{eq:key2}, thanks to \eqref{eq:DamGLw}, we get
	\begin{equation}\label{eq:lhs}
		\begin{split}
			&
			\int |\grad u^-|^2 \left[ 1 + \frac{|\lambda_2/G_d|}{\sqrt{\lambda_1/G_d}} \kor \right] 
			-
			\frac{Q-1}{2}
			\frac{|\lambda_2/G_d|}{\sqrt{\lambda_1/G_d}}
			\int \gl^2 \frac{|u^-|^2}{\kor}
			\\
			&\geq\,
			\int |\grad u^-|^2 
			+
			\frac{Q-3}{Q-1}\frac{|\l_2/G_d|}{\sqrt{\l_1/G_d}}
			\int 
			\kor
			|\grad u^-|^2
			\\
			&\geq\,
			\norm{\grad u^-}^2
			.
		\end{split}
	\end{equation}
	
	In order to estimate the right-hand side of \eqref{eq:key2}, first of all observe that choosing $\Psi_2\equiv 1$ in \eqref{eq:fondws2}, and using Lemma\til\ref{lem:B} with $\gamma=1$, we have
	\begin{equation*}
		\frac{\l_2}{G_d} \int |\mathcal{J}_\phi| \frac{1}{\jap{s}}|\u|^2 
		= 
		\l_2 \int |\mathcal{J}_\phi| |\u|^2 
		= 
		\im \int |\mathcal{J}_\phi| \f \o{\u}
	\end{equation*}
	which gives, changing back the variables,
	\begin{equation*}
		\frac{\l_2}{G_d} \int \gl^2 |u|^2 
		= 
		\im \int f \o{u}
		.
	\end{equation*}
	Therefore we may estimate, by H\"older's and Young's inequalities,
	\begin{equation}\label{eq:optimizing-1}
			\sqrt{\frac{|\l_2|}{G_d}} \norm{\gl u}
			\le 
			\sqrt{\int |f u|}
			\le
			\norm{ \frac{\kor}{\gl} f}^{\frac12} \norm{ \gl \frac{u}{\kor} }^{\frac12}
			\leq
			\frac1{2\gamma} \norm{\frac{\kor}{\gl} f}
			+
			\frac\gamma2 \norm{ \gl \frac{u}{\kor} }
	\end{equation}
	for any $\gamma>0$.	
	By \eqref{eq:optimizing-1} and the Hardy inequality \eqref{eq:hardyGL}, since $|\l_2| \le \delta \l_1$, we obtain
	\begin{equation}\label{eq:rhs}
		\begin{split}
			&
			\left|
				\re \int \frac{\kor}{\gl} \o{e^{-i\psi} f} 
				\left[ (Q-1) \gl \frac{u^-}{\kor}
				+ 
				\frac{|\lambda_2/G_d|}{\sqrt{\lambda_1/G_d}} \gl u^-
				+ 
				2 \left[ \frac{|z|^2}{\kor^2} \hat{z} + \frac{t}{\kor^2} \hat{z}^\perp \right] \cdot \grad u^- \right]
			\right|
			\\
			&\le\,
			\norm{\frac{\kor}{\gl} f}
			\left[
				(Q-1) \norm{\gl \frac{u^-}{\kor}} 
				+
				\frac{\sqrt{\delta}}{2\gamma} \norm{\frac{\kor}{\gl} f}
				+
				\frac{\sqrt{\delta}\gamma}{2} \norm{\gl \frac{u^-}{\kor}}
				+
				2 \norm{\grad u^-}
			\right]
			\\
			&\le\,
			\norm{\frac{\kor}{\gl} f}
			\left[
			\frac{2(Q-1) + \sqrt{\delta}\gamma + 2(Q-2)}{Q-2} 
			\norm{ \grad u^- }
			+
			\frac{\sqrt{\delta}}{2\gamma} \norm{\frac{\kor}{\gl} f}
			\right]
			\\
			&=\,
			\frac{8d+2+\gamma\sqrt{\delta}}{2d}
			\norm{\frac{\kor}{\gl} f} \norm{\grad u^-}
			+
			\frac{\sqrt{\delta}}{2\gamma}
			\norm{\frac{\kor}{\gl} f}^2
		\end{split}
	\end{equation}
	for any $\gamma>0$.
	From the above estimate, \eqref{eq:key2} and \eqref{eq:lhs}, we finally obtain
    \begin{equation}\label{eq:est-final}
         \norm{\grad u^-}^2-\frac{8d+2+\gamma\sqrt{\delta}}{2d}
		\norm{ \frac{\kor}{\gl} f} 
		\norm{\grad u^-}
		-\frac{\sqrt{\delta}}{2\gamma}
		\norm{ \frac{\kor}{\gl} f}^{2}
		\le
		0.
	\end{equation}
	Observe that the quadratic inequality above represents a parabola in the $ \norm{\grad u^-}$ variable, hence it implies that
	\begin{equation*}
	\norm{\grad u^-}
	\leq 
	\left(
		\frac{8d+2+\gamma\sqrt{\delta}}{4d}+\sqrt{\frac{(8d+2+\gamma\sqrt{\delta})^2}{16d^2}+\frac{\sqrt{\delta}}{2\gamma}} 
	\right) 
	\norm{ \frac{\kor}{\gl} f}
	\end{equation*}
	for any $\gamma>0$, from which \eqref{eq:est2} follows.

\subsubsection{Removing regularity.}
So far we simply assumed that $u$ and $f$ are smooth and with compact support. Now we remove this additional assumption, applying a standard cut-off and mollification argument. 
	Consider a smooth cut-off function $\chi \colon \R_+ = [0,+\infty) \to [0,1]$ such that
	\begin{equation}\label{eq:prop-chi}
		\begin{split}
			\chi(r) \equiv 1 \quad\text{for $r\le 1$},&
			\qquad
			\chi(r) \equiv 0 \quad\text{for $r \ge 2$},
			\\
			\chi' \in L^\infty(\R_+),&
			\qquad
			r\chi'' \in L^\infty(\R_+).
		\end{split}
	\end{equation}
	Define the scaled cut-off function, for $R>0$,
	\begin{equation*}
		\chi_R(z,t) := \chi\left( \frac{|(z,t)|_{\H}}{R} \right).
	\end{equation*}	
Recalling \eqref{eq:kder1}, we have 
	\begin{align*}
		\grad \, \chi_R(z,t) &=
		\frac{1}{R}
		\chi'\left( \frac{|(z,t)|_{\H}}{R} \right)
		\frac{|z|^2 z + (y,-x) t}{|(z,t)|_{\H}^{3}}
		\\
		\L \, \chi_R(z,t) &=
		-
		\frac{1}{R}
		\left[
		\frac{|(z,t)|_{\H}}{R} \chi''\Big(\frac{|(z,t)|_{\H}}{R} \Big) + (2d+1) \chi'\Big(\frac{|(z,t)|_{\H}}{R}\Big) \right]
		\frac{|z|^2}{|(z,t)|_{\H}^3}
	\end{align*}
	and hence, due to properties \eqref{eq:prop-chi} and to \eqref{eq:kder2}, it follows
	\begin{equation}\label{eq:prop-chiR}
		\begin{split}
			\chi_R(z,t) \equiv 1 \quad\text{for $|(z,t)|_{\H}\le R$},&
			\qquad
			\chi_R(z,t) \equiv 0 \quad\text{for $|(z,t)|_{\H} \ge 2R$},
			\\
			|\grad \, \chi_R(z,t)| 
			\le
			C_{\chi,d}
			\frac{1}{R} \gl
			,
			&
			\qquad
			|\L \, \chi_R(z,t)|
			\le
			C_{\chi,d}
			\frac{1}{R}
			\frac{\gl}{\kor}
			,
		\end{split}
	\end{equation}
	where $C_{\chi,d}>0$ is a constant depending only on $\chi$ and $d$.
	Let now $\varphi\in\mathcal C^{\infty}_0([0,+\infty);[0,+\infty))$ be a positive decreasing function such that the kernel on $\H^d$ defined by $\phi(z,t)=\varphi(|(z,t)|_{\H})$ has the property $\int_{\H^d}\phi(z,t)\,dz\,dt=1$. By scaling, we construct the delta-sequence
	$$\phi_\e(z,t) := \e^{-2d-1} \phi\left(\frac{z}{\e},\frac{t}{\e^2}\right)$$
	such that $\lim_{\e\to0^+} \phi_\e(z,t) = \delta(z,t)$, being $\delta$ the Dirac mass centered at the origin. If $u\in H^1(\H^d)$ is a weak solution to \eqref{eq:solution}, then 
	\begin{equation*}
		u_{R,\e} := \chi_R u * \phi_\e
	\end{equation*}
	solves weakly the equation
	\begin{equation*}
		- \L u_{R,\e} + \lambda u_{R,\e} =
		\left(\chi_R f - 2 \grad \chi_R \cdot \grad u + u \L \chi_R \right)*\phi_\e
		=: \widetilde{f}_{R,\e},
	\end{equation*}
	namely
	\begin{equation*}
		- \langle \grad v, \grad u_{R,\e} \rangle
		+
		\lambda
		\langle v, u_{R,\e} \rangle
		=
		\langle v, \widetilde{f}_{R,\e} \rangle
	\end{equation*}
	for any $v\in H^1(\H^d)$.
	From the definition of $\widetilde{f}_{R,\e}$ and the properties \eqref{eq:prop-chiR}, note that
	\begin{align}\label{eq:above}
		\norm{\frac{\kor}{\gl} \widetilde{f}_{R,\e}}
		&\le
		\norm{ \frac{\kor}{\gl} \left[ \chi_R f - 2 \grad \chi_R \cdot \grad u + u \L \chi_R \right]}
		\\
		&\le
		\norm{ \frac{\kor}{\gl} f}
		+
		4 C_{\chi,d} \left(\int_{R<\kor<2R} |\grad u|^2 \right)^{\frac12}
		+
		\frac{C_{\chi,d}}{R} 	
		\left(\int_{R<\kor<2R} |u|^2 \right)^{\frac12}
		\nonumber
		\\
		&=
		\norm{\frac{\kor}{\gl} f} + o(1),
		\nonumber
	\end{align}
	when $R\to+\infty$, since $u\in H^1(\H^d)$.
	Finally, since $u_{R,\e}$ and $\widetilde{f}_{R,\e}$ are smooth and compactly supported, \eqref{eq:est-final} holds true with $\widetilde{f}_{R,\e}$ and $u_{R,\e}^-(z,t) := e^{-i\sgn(\l_2)\sqrt{\l_1/G_d} \Kor}u_{R,\e}(z,t)$ in place of $f$ and $u^-$ respectively. By \eqref{eq:above}, we finally get
	\begin{equation*}
		\big\|\grad u_{R,\e}^-\big\|
		\le
		\left(
			\frac{8d+2+\gamma\sqrt{\delta}}{4d}+\sqrt{\frac{(8d+2+\gamma\sqrt{\delta})^2}{16d^2}+\frac{\sqrt{\delta}}{2\gamma}} 
		\right) 
		\left( \norm{ \frac{\kor}{\gl} f} +o(1) \right),
	\end{equation*}
	when $R\to+\infty$.
	In conclusion, letting first $\e\to0^+$ for fixed $R>0$ and then $R \to +\infty$, by the dominated convergence and monotonic convergence theorems we complete the proof of \eqref{eq:est2}.
	
\section{Proof of Theorem \ref{thm:V1}}\label{sec:proof2}

We are now ready to prove Theorem \ref{thm:V1}. We need to prove that the equation
\begin{equation}\label{eq:prob1}
- \L u-Vu+\lambda u=0
\end{equation}
has no non-trivial radial solutions in $H^1(\H^d)$, for any $\lambda\in\C$, so that we can consider \eqref{eq:prob0} with $f=Vu$.
We will treat separately the two cases $|\lambda_2|>\delta_* \lambda_1$, and $|\lambda_2|\leq \delta_* \lambda_1$, where $\delta_*>0$ is implicitly defined in \eqref{eq:delta*}. 

\subsection{Case $|\lambda_2|> \delta_* \lambda_1$}
When $\lambda_2>\delta_*\lambda_1$, arguing as in \eqref{eq:|l2|>l1-12} we obtain
\begin{equation*}\label{eq:|l2|>l1-123}
	0 \ge (\delta_* \l_1-\l_2) \norm{u} \ge 
	\delta_* \norm{\grad u} \left(\norm{\grad u} - \frac{1+1/\delta_*}{d} \norm{ \frac{\kor}{\gl} Vu} \right)
\end{equation*}
and by the assumption \eqref{eq:assV1}  we infer that
\begin{equation*}
	0
	\geq
	\left(
		1 - b  \frac{1+1/\delta_*}{d}
	\right)
	\norm{\grad u}.
\end{equation*}
By the identity \eqref{eq:kd-def2} it is easy to realize that \eqref{eq:nuovanuova} is equivalent to $1-b\frac{1+1/\delta_*}{d}>0$, therefore we necessarily conclude that $u\equiv0$. The case $\lambda_2<-\delta_* \lambda_1$ is completely analogous.
\subsection{Case $|\lambda_2|\leq \delta_* \lambda_1$}
As above, we read \eqref{eq:prob1} as \eqref{eq:prob0} with $f=Vu$. We then go back to the key identity \eqref{eq:key2}, which now reads
\begin{multline*}
 \int |\grad u^-|^2 \left[ 1 + \frac{|\lambda_2/G_d|}{\sqrt{\lambda_1/G_d}} \kor \right] 
		-
		\frac{Q-1}{2}
		\frac{|\lambda_2/G_d|}{\sqrt{\lambda_1/G_d}}
		\int \gl^2 \frac{|u^-|^2}{\kor}
		\\
		=
		-
		\re \int \frac{\kor}{\gl} \o{Vu^-} 
		\left[ (Q-1) \gl \frac{u^-}{\kor}
		+ 
		\frac{|\lambda_2/G_d|}{\sqrt{\lambda_1/G_d}} \gl u^-
		+ 
		2 \left[ \frac{|z|^2}{\kor^2} \hat{z} + \frac{t}{\kor^2} \hat{z}^\perp \right] \cdot \grad u^- \right]
		.
\end{multline*}
	The left-hand side is estimated as in \eqref{eq:lhs}. For the right-hand side, we argue in a slightly different way than in \eqref{eq:optimizing-1} and we involve assumption \eqref{eq:assV1}, which together with the Hardy inequality \eqref{eq:hardyGL} gives 
	\begin{equation*}
		\sqrt{\frac{|\lambda_2|}{G_d}}
		\norm{ \gl u} \le \sqrt{\int |V||u|^2}
		\le
		\norm{ \frac{\kor}{\gl} Vu^-}^{1/2} \norm{\gl \frac{u^-}{\kor}}^{1/2}\leq \sqrt{\frac{b}{d}} \norm{\grad u^-}.
	\end{equation*}
	By the above inequality, assumption \eqref{eq:assV1}, and the Hardy inequality \eqref{eq:hardyGL}, since $|\l_2| \le \delta_* \l_1$ we can estimate
	\begin{equation}\label{eq:rhs2}
		\begin{split}
		&
		\left|
			\re \int \frac{\kor}{\gl} \o{Vu^-} 
			\left[ (Q-1) \gl \frac{u^-}{\kor}
			+ 
			\frac{|\lambda_2/G_d|}{\sqrt{\lambda_1/G_d}} \gl u^-
			+ 
			2 \left[ \frac{|z|^2}{\kor^2} \hat{z} + \frac{t}{\kor^2} \hat{z}^\perp \right] \cdot \grad u^- \right]
		\right|
		\\
		&\le\,
		\norm{\frac{\kor}{\gl} Vu^-}
		\left[
			(Q-1) \norm{\gl \frac{u^-}{\kor}}
			+
			\frac{|\lambda_2/G_d|}{\sqrt{\lambda_1/G_d}} \norm{\gl u^-}
			+
			2 \norm{\grad u^-}
		\right]
		\\
		&\le\,
		\left( b^{3/2} \sqrt{\frac{\delta_*}{d}} + b \frac{4d+1}{d} \right)
		\norm{\grad u^-}^2.
	\end{split}
	\end{equation}
	By \eqref{eq:lhs} and \eqref{eq:rhs2}, we get
	\begin{equation}\label{eq:nuova}
	\left( 1 - b^{3/2} \sqrt{\frac{\delta_*}{d}} - b \frac{4d+1}{d} \right)
	\norm{\grad u^-}^2+
	\frac{2d-1}{2d+1}\frac{|\l_2/G_d|}{\sqrt{\l_1/G_d}}
		\int 
		\kor
		|\grad u^-|^2
	\leq 0.
	\end{equation}
	Since from \eqref{eq:Kd-impdef} and \eqref{eq:kd-def2} it can be seen that
	\begin{equation*}
		\left[ \frac{1}{d \kappa_d} \right]^{3/2} \sqrt{\frac{\delta_*}{d}} + \left[ \frac{1}{d\kappa_d} \right] \frac{4d+1}{d} = 1,
	\end{equation*}
	it is easy to conclude that  \eqref{eq:nuovanuova} implies that the coefficient of $\|\nabla_{\!\H} u^-\|^2$ in \eqref{eq:nuova} is positive, so that necessarily $u\equiv0$. This completes the proof of Theorem \ref{thm:V1}.
	%
	
\section{Proof of Theorems \ref{thm:pp} and \ref{thm:grad-est2}}\label{sec:proof3}

We now prove Theorems \ref{thm:pp} and \ref{thm:grad-est2}, the proof of which is analogous to the one of Theorem \ref{thm:grad-est}, taking into account the presence of $V$. 

\subsection{Proof of estimates \eqref{eq:est3b} and \eqref{eq:est3}} 
\label{subse:proofest3}
As above, let us first consider the case $\lambda_2> \delta \lambda_1$. Choosing $v=u$ in \eqref{eq:solution2} as test function, subtracting $\delta$ times the real part of the resulting identity from its imaginary part, we get
	\begin{equation}\label{eq:|l2|>l1-100}
		(\delta \l_1-\l_2) \int |u|^2 = \delta \int |\grad u|^2 + \delta \int V|u|^2 + \delta \re \int \o{u} f - \im \int \o{u} f.
	\end{equation}
	By Cauchy--Schwarz and the Hardy inequality \eqref{eq:hardyGL}, we estimate
	\begin{equation*}
		\Big\lvert 
		\delta
		\re \int \o{u} f - \im \int \o{u} f
		\Big\rvert
		\le
		(1+\delta) \int |u||f|
		\le
		(1+\delta) \norm{ \frac{\kor}{\gl} f } \norm{ \gl \frac{u}{\kor} }
		\le
		\frac{1+\delta}{d} \norm{ \frac{\kor}{\gl} f} \norm{\grad u}.
	\end{equation*}	
	Plugging the last inequality in \eqref{eq:|l2|>l1-100}, neglecting the positive term involving $V_+$ and using the subordination condition \eqref{eq:b1}, we get
	\begin{equation*}\label{eq:|l2|>l1-124}
	0 \ge (\delta \l_1-\l_2) \norm{u}^2 \ge 
	\delta
	\norm{\grad u} \left( (1-b_1^2)\norm{\grad u} - \frac{1+1/\delta}{d} \norm{ \frac{\kor}{\gl} f} \right)
	\end{equation*}
(in the case of the estimate \eqref{eq:est3b}, we have $b_1=0$) from which, if $u$ is non-vanishing, estimates \eqref{eq:est3b} and \eqref{eq:est3} necessarily follow. The case $\lambda_2<-\delta \lambda_1$ is completely analogous and we omit the details.
\subsection{Proof of estimate \eqref{eq:est4b}}

Arguing as in the proof of \eqref{eq:est2}, and involving $V$ in the algebraic manipulations, after some integration by parts one obtains the analogous to \eqref{eq:key2}, which is (recall the definition of $G_d$ in \eqref{eq:Gd})
\begin{equation}\label{eq:4b1}
	\begin{split}
		& \int \left[ 1 + \frac{|\lambda_2/G_d|}{\sqrt{\lambda_1/G_d}} \kor \right] |\grad u^-|^2
		-
		\frac{Q-1}{2}
		\frac{|\lambda_2/G_d|}{\sqrt{\lambda_1/G_d}}
		\int \gl^2 \frac{|u^-|^2}{\kor}
		\\
		&=
		-
		\re \int \frac{\kor}{\gl} \o{e^{-i\psi}(f+Vu)} 
		\left\{
			\left[
				\frac{Q-1}{\kor} 
				+ 
				\frac{|\lambda_2/G_d|}{\sqrt{\lambda_1/G_d}} 
			\right]
			\gl u^-
			+ 
			2 \left[ \frac{|z|^2}{\kor^2} \hat{z} + \frac{t}{\kor^2} \hat{z}^\perp \right] \cdot \grad u^- 
		\right\}
		.
	\end{split}
\end{equation}
Now, notice that, for a generic function $g$,
\begin{equation*}
	\left[ \frac{|z|^2}{\kor^2} \hat{z} + \frac{t}{\kor^2} \hat{z}^\perp \right] \cdot \grad g
	=
	\frac{\grad \kor}{|\grad \kor|} \cdot \grad g
	=
	\frac{\gl}{\kor} \left[ z \cdot \nabla_{z} + 2t \partial_t \right] g + \frac{t}{\kor^2} \hat{z}^\perp \cdot \nabla_z g
	=
	\gl \partial_\kor g + \frac{t}{\kor^2} \hat{z}^\perp \cdot \nabla_z g
	.
\end{equation*}
If in particular $g(z,t) = g^0(\kor)$ is a radial function, then
\begin{equation*}
	\left[ \frac{|z|^2}{\kor^2} \hat{z} + \frac{t}{\kor^2} \hat{z}^\perp \right] \cdot \grad g
	=
	\gl \partial_\kor g 
	.
\end{equation*}
Therefore the identity \eqref{eq:4b1} becomes
\begin{equation}\label{eq:keyV}
	\begin{split}
		& \int \left[ 1 + \frac{|\lambda_2/G_d|}{\sqrt{\lambda_1/G_d}} \kor \right] |\grad u^-|^2
		-
		\frac{Q-1}{2}
		\frac{|\lambda_2/G_d|}{\sqrt{\lambda_1/G_d}}
		\int \gl^2 \frac{|u^-|^2}{\kor}
		\\
		&
		+ 
		\frac{|\lambda_2/G_d|}{\sqrt{\lambda_1/G_d}} \int \kor V |u^-|^2
		-
		\int ( V + \kor \partial_{\kor} V ) |u^-|^2
		\\
		&=
		-
		\re \int \frac{\kor}{\gl} \o{e^{-i\psi} f} 
		\left\{
		\left[
		\frac{Q-1}{\kor} 
		+ 
		\frac{|\lambda_2/G_d|}{\sqrt{\lambda_1/G_d}} 
		\right]
		\gl u^-
		+ 
		2 \left[ \frac{|z|^2}{\kor^2} \hat{z} + \frac{t}{\kor^2} \hat{z}^\perp \right] \cdot \grad u^- 
		\right\}
		.
	\end{split}
\end{equation}
%
%
	For the estimate of the right-hand side of \eqref{eq:keyV} we proceed exactly as in the free case $V\equiv0$, and we have \eqref{eq:rhs}.
	In order to estimate the left-hand side we need some slight modifications, due to the presence of $V$.
	First of all, by \eqref{eq:b2b} we obtain
	\begin{equation}\label{eq:iuno}
	-\int\left(V+\kor \partial_{\kor} V \right)|u^-|^2
	=
	-\int \partial_{\kor} \left( \kor V \right)|u^-|^2
	\geq
	-\int \left[ \partial_{\kor} (\kor V) \right]_+ |u^-|^2
	\geq-b^2\int|\nabla_{\!\H}u^-|^2.
	\end{equation}
	Then, arguing as above, we arrive at the estimate  
	 \begin{equation*}
		(1-b^2)\norm{\grad u^-}^2
		-\frac{8d+2+\gamma\sqrt{\delta}}{2d}
		\norm{ \frac{\kor}{\gl} f} 
		\norm{\grad u^-}
		-\frac{\sqrt{\delta}}{2\gamma}
		\norm{ \frac{\kor}{\gl} f}^{2}
		\le
		0.
	\end{equation*}
The proof now follows exactly the same as that of Theorem \ref{thm:grad-est}.

\subsection{Proof of estimate \eqref{eq:est4}}
Let us start again by the key identity \eqref{eq:keyV}. 
	For the estimate of the right-hand side, we have \eqref{eq:rhs} as above.  We now only need to estimate the additional term containing $V$, for which, involving  \eqref{eq:b1}, we get 
	\begin{equation}\label{eq:idue0}
		\int \kor V|u^-|^2 
		\geq
		-
		\int V_-\big| \sqrt{\kor} \, u^-\big|^2
		\geq
		-b_1^2
		\int\left|\nabla_{\!\H}\left( \sqrt{\kor} \, u^-\right)\right|^2.
	\end{equation}
	Moreover, recalling \eqref{eq:kder2}, the following identity holds: 
	\begin{equation}\label{eq:idue1}
		\int \kor |\grad u^-|^2
		- \frac{Q-1}{2} \int \frac{\gl^2}{\kor} |u^-|^2
		=
		\int |\grad(\sqrt{\kor} \, u^-)|^2
		-
		\frac14 \int \frac{\gl^2}{\kor} |u^-|^2
		.
	\end{equation}
	Putting together \eqref{eq:iuno}, \eqref{eq:idue0} and \eqref{eq:idue1} we get
	\begin{equation}\label{eq:lhsquasi}
	\begin{split}
		& \int \left[ 1 + \frac{|\lambda_2/G_d|}{\sqrt{\lambda_1/G_d}} \kor \right] |\grad u^-|^2
		-
		\frac{Q-1}{2}
		\frac{|\lambda_2/G_d|}{\sqrt{\lambda_1/G_d}}
		\int \gl^2 \frac{|u^-|^2}{\kor}
		\\
		&
		+ 
		\frac{|\lambda_2/G_d|}{\sqrt{\lambda_1/G_d}} \int \kor V |u^-|^2
		-
		\int ( V + \kor \partial_{\kor} V ) |u^-|^2
		\\
		&\geq\,
		(1-b_2^2)\int|\grad u^-|^2
		+
		(1-b_1^2) \frac{|\l_2/G_d|}{\sqrt{\l_1/G_d}} 
		\int | \grad( \sqrt{\kor} \, u^- ) |^2
		-
		\frac14 \frac{|\l_2/G_d|}{\sqrt{\l_1/G_d}} \int \gl^2 \frac{|u^-|^2}{\kor}.
	\end{split}
	\end{equation} 
To complete the estimate of the left-hand side of \eqref{eq:keyV} and prove its positivity, we finally need to control the last negative term in \eqref{eq:lhsquasi}. First, notice that \eqref{eq:fond2} holds also in this case, since $V\in\R$, and arguing as in \eqref{eq:optimizing-1} we have that
\begin{equation}\label{eq:optimizingV}
		\sqrt{\frac{|\l_2|}{G_d}} \norm{\gl u}
		\leq 
		\frac1{2\widetilde{\gamma}} \norm{\frac{\kor}{\gl} f} 
		+
		\frac{\widetilde{\gamma}}{2} \norm{\frac{\gl}{\kor} u}
	\end{equation}
	for any $\widetilde{\gamma}>0$. By this and the Hardy inequality \eqref{eq:hardyGL}, and since $|\lambda_2|\leq \delta \lambda_1$, it follows
\begin{equation*}\label{eq:itre}
	\frac{|\l_2/G_d|}{\sqrt{\l_1/G_d}} \int \gl^2 \frac{|u^-|^2}{\kor}
	\leq
	\frac{|\l_2/G_d|}{\sqrt{\l_1/G_d}} \norm{ \gl u^- } \norm{\gl \frac{u^-}{\kor}}
	\leq
	\frac{\sqrt{\delta}}{2\widetilde{\gamma}d}
	\norm{\frac{\kor}{\gl} f} \norm{\grad u^-}
	+
	\frac{\widetilde{\gamma}\sqrt{\delta}}{2d^2} \norm{\grad u^-}^2
	.
\end{equation*}	
Plugging this information in \eqref{eq:lhsquasi} we obtain
\begin{equation}\label{eq:lhsquasigia}
	\begin{split}
		& \int \left[ 1 + \frac{|\lambda_2/G_d|}{\sqrt{\lambda_1/G_d}} \kor \right] |\grad u^-|^2
		-
		\frac{Q-1}{2}
		\frac{|\lambda_2/G_d|}{\sqrt{\lambda_1/G_d}}
		\int \gl^2 \frac{|u^-|^2}{\kor}
		\\
		&
		+ 
		\frac{|\lambda_2/G_d|}{\sqrt{\lambda_1/G_d}} \int \kor V |u^-|^2
		-
		\int ( V + \kor \partial_{\kor} V ) |u^-|^2
		\\
		&\geq\,
		\left( 1-b_2^2 - \frac{\widetilde{\gamma}\sqrt{\delta}}{8d^2}\right) \norm{\grad u^-}^2
		-
		\frac{\sqrt{\delta}}{8d \widetilde{\gamma}} \norm{\frac{\kor}{\gl} f} \norm{\grad u^-}
		\\
		&+
		(1-b_1^2) \frac{|\l_2/G_d|}{\sqrt{\l_1/G_d}} 
		\int | \grad( \sqrt{\kor} \, u^- ) |^2
		.
	\end{split}
\end{equation} 
	Finally, by \eqref{eq:rhs} and \eqref{eq:lhsquasigia} we conclude, since $0\le b_1<1$ by assumption, that
	\begin{equation*}\label{eq:cisiamo}
	\left(1-b_2^2-\frac{\widetilde{\gamma}\sqrt{\delta}}{8d^2} \right)
	\norm{\grad u^-}^2
	-
	\left(\frac{8d+2+\gamma\sqrt{\delta}}{2d}+\frac{\sqrt{\delta}}{8d\widetilde{\gamma}}\right)
	\norm{\frac{\kor}{\gl} f} \norm{\grad u^-}
	-
	\frac{\sqrt{\delta}}{2\gamma} \norm{\frac{\kor}{\gl} f}^2
	\leq
	0
	\end{equation*}
	for any $\gamma,\widetilde{\gamma}>0$. 
	Since we need the coefficient of $\norm{\grad u^-}^2$ to be positive, we impose also $\widetilde{\gamma} < 8d^2 (1-b_2^2)/\sqrt{\delta}$ (remember that by assumption $0\le b_2<1$).
	After a change of variables $\gamma_2^2 := \frac{\widetilde{\gamma}\sqrt{\delta}}{8d^2(1-b_2^2)}$ and renaming $\gamma_1 := \gamma$, the last inequality above implies
	$$
	\norm{\grad u^-}
	\leq
	g_{d,\delta,b_2}(\gamma_1,\gamma_2) \norm{\frac{\kor}{\gl} f},
	$$
	for any $\gamma_1>0$ and $0<\gamma_2 < 1$,
	where $g_{d,\delta,b_2}(\gamma_1,\gamma_2)$ is defined in the statement of Theorem\til\ref{thm:grad-est2}.
	
	This formally completes the proof of estimate \eqref{eq:est4}.
	To make the argument rigorous, one can argue as in the free case. Nevertheless, some additional regularity needs to be assumed on $V$ in order to close the regularization argument, that is condition $V\in W_{\operatorname{loc}}^{1,p}(\H^d)$ in the statement of Theorem \ref{thm:grad-est2}. We address the reader to the proof of \cite[Theorem 3.1]{CFK}, in which this topic is treated in detail in the Euclidean case, for which the regularization works in a completely analogous way as in the case of the present manuscript.
\section{Proof of Theorem \ref{thm:V2}}\label{sec:proof4}
The proof of Theorem \ref{thm:V2} follows exactly the same line as the proofs of the rest of the theorems, so we only sketch it. First of all, 
we write the resolvent equation as
\begin{equation}
- \L u-(\re V)u+\lambda u=i\im V u,
\end{equation}
and we aim to prove that the only $H^1$-solution is the null one. Let us fix $\widetilde{\delta} > 0$ such that
\begin{equation}\label{eq:delta**}
	\sqrt{ \frac{b_3}{d(1-b_1^2)}}
	<
	\sqrt{\widetilde{\delta}}
	<
	\sqrt{\frac{d}{b_3}} \left(1-b_2^2-\frac{4d+1}{d}b_3 \right) \cdot
	\left(\frac{1}{4d}+b_3\right)^{-1}
	.
\end{equation}
This choice is possible thanks to our assumption \eqref{eq:assf.Lambdanew}. Indeed, observe that the inequality
\begin{equation*}
	\sqrt{ \frac{b_3}{d(1-b_1^2)}}
	<
	\sqrt{\frac{d}{b_3}} \left(1-b_2^2-\frac{4d+1}{d}b_3 \right) \cdot
	\left(\frac{1}{4d}+b_3\right)^{-1}
\end{equation*}
is equivalent to
\begin{equation*}
	b_3^2 + \left( \frac{1}{4d} + (4d+1) \sqrt{1-b_1^2} \right) b_3 - d (1-b_2^2) \sqrt{1-b_1^2} < 0
\end{equation*}
and the left-hand side represents a parabola in the variable $b_3$, so assumption \eqref{eq:assf.Lambdanew} arises after computing the positive root.

The easier case is when $|\lambda_2|>\widetilde{\delta}\lambda_1$, for which we can argue like in the analogous case in Subsection\til\ref{subsec:est1}. Indeed, take $f=Vu$, by Cauchy--Schwarz and Hardy inequalities, and by the subordination conditions \eqref{eq:assV12} and \eqref{eq:useful}, we have
	\begin{equation*}
		\widetilde{\delta} \re \int |u|^2V - \im \int |u|^2 V
		\ge
		-
		\widetilde{\delta} \Big(b_1^2+b_3\cdot\frac{1/{\widetilde{\delta}}}{d} \Big) \int |\grad{u}|^2.
	\end{equation*}	
Analogously as in Subsection \ref{subsec:est1}, in view of \eqref{eq:|l2|>l1-1} and of the latter inequality, we get
	\begin{equation*}\label{eq:|l2|>l1-125}
	0 \ge (\widetilde{\delta}\l_1-\l_2) \norm{u}^2 
	\ge 
	\widetilde{\delta} \norm{\grad u}\Big( 1 - b_1^2 - b_3 \cdot \frac{1/\widetilde{\delta}}{d}\Big)
	\end{equation*}
	from which, if $u$ is non-vanishing, then necessarily
	\begin{equation*}
		\Big( 1 - b_1^2 - b_3 \cdot \frac{1/\widetilde{\delta}}{d}\Big) \norm{\grad u} \le 0.	
	\end{equation*}
Notice that
$$
 1 - b_1^2 - b_3 \cdot \frac{1/\widetilde{\delta}}{d}
\le 0 \Longleftrightarrow \frac{b_3}{d(1-b_1^2)}
	\ge
	\widetilde{\delta},
$$
hence in view of the left-hand side of \eqref{eq:delta**}, we conclude that $u\equiv 0$ and hence the thesis of the theorem in this case.

In the cone $|\lambda_2|\leq\widetilde{\delta}\lambda_1$, we start with the identity \eqref{eq:keyV} replacing $V$ with $\re V$ and setting $f=i(\im V)u$, which reads
\begin{equation}\label{eq:keyV2}
	\begin{split}
		& \int \left[ 1 + \frac{|\lambda_2/G_d|}{\sqrt{\lambda_1/G_d}} \kor \right] |\grad u^-|^2
		-
		\frac{Q-1}{2}
		\frac{|\lambda_2/G_d|}{\sqrt{\lambda_1/G_d}}
		\int \gl^2 \frac{|u^-|^2}{\kor}
		\\
		&
		+ 
		\frac{|\lambda_2/G_d|}{\sqrt{\lambda_1/G_d}} \int \kor \re V |u^-|^2
		-
		\int ( \re V + \kor \partial_{\kor} \re V ) |u^-|^2
		\\
		&=
		-
		\im \int \frac{\kor}{\gl} \im V 
		\left\{
		\left[
		\frac{Q-1}{\kor} 
		+ 
		\frac{|\lambda_2/G_d|}{\sqrt{\lambda_1/G_d}} 
		\right]
		\gl |u^-|^2
		+ 
		2 \left[ \frac{|z|^2}{\kor^2} \hat{z} + \frac{t}{\kor^2} \hat{z}^\perp \right] \cdot \o{u^-} \grad u^- 
		\right\}
		.
	\end{split}
\end{equation}
	For the left-hand term, we proceed almost exactly as in the proof of estimate \eqref{eq:est4}, with the only difference that in place of \eqref{eq:optimizingV} we use instead
	\begin{equation}\label{eq:optimizingV2}
		\sqrt{\frac{|\l_2|}{G_d}}
		\norm{\gl u} \le \sqrt{\int |\im V| |u|^2}
		\le
		\norm{ \frac{\kor}{\gl} \im V u}^{1/2} \norm{\gl \frac{u}{\kor}}^{1/2}\leq \sqrt{\frac{b_3}{d}}\|\nabla_{\!\H} u^-\|
	\end{equation}
	obtained combining \eqref{eq:fond2}, \eqref{eq:assV2} and the Hardy inequality \eqref{eq:hardyGL}.
Consequently 
\begin{equation}\label{eq:iquattro}
\frac{|\l_2/G_d|}{\sqrt{\l_1/G_d}}
\int \gl^2 \frac{|u^-|^2}{\kor}
\leq
\frac{|\l_2/G_d|}{\sqrt{\l_1/G_d}} \norm{\gl u^-} \norm{\gl \frac{u^-}{\kor}}
\leq
\frac{\sqrt{\widetilde{\delta} b_3}}{d^{3/2}}\|\nabla_{\!\H} u^-\|^2,
\end{equation}	
since $|\lambda_2|\leq \widetilde{\delta}\lambda_1$. 
For the right-hand side, we proceed as in the proof of Theorem\til\ref{thm:V1}, obtaining, thanks to the inequalities \eqref{eq:iquattro} and \eqref{eq:useful}, an estimate analogous to \eqref{eq:rhs2}, but with $V$, $b$ and $\delta_*$ substituted by $\im V$, $b_3$ and $\widetilde{\delta}$ respectively.
Putting together all the informations we finally obtain the estimate
\begin{equation*}
	\left( 
		1 - b_2^2 - b_3^{3/2} \cdot \sqrt{ \frac{\widetilde{\delta}}{d} } - b_3 \cdot \frac{4d+1}{d} - b_3^{1/2} \cdot \frac{1}{4d} \sqrt{ \frac{\widetilde{\delta}}{d} }
	\right)
	\| \grad u^-\|^2
	+(1-b_1^2)\frac{|\l_2/G_d|}{\sqrt{\l_1/G_d}}
	\int
	\big| \grad \big( \sqrt{\kor} \, u^-\big)\big|^2\leq0,
\end{equation*}	
and in particular, since $b_1<1$,
\begin{equation*}\label{eq:ultima}
	\left( 
	1 - b_2^2 - b_3^{3/2} \cdot \sqrt{ \frac{\widetilde{\delta}}{d} } - b_3 \cdot \frac{4d+1}{d} - b_3^{1/2} \cdot \frac{1}{4d} \sqrt{ \frac{\widetilde{\delta}}{d} }
	\right)
	\norm{\grad u^-}^2
\leq0.
\end{equation*}
The thesis now follows from \eqref{eq:delta**}.


\section*{Acknowledgement}

L. Fanelli, L. Roncal and N. M. Schiavone 
are partially supported 
by the Basque Government through the BERC 2022--2025 program 
and 
by the Spanish Agencia Estatal de Investigación
through BCAM Severo Ochoa excellence accreditation CEX2021-001142-S/MCIN/AEI/10.13039/501100011033.

L. Fanelli is also supported by the projects PID2021-123034NB-I00/MCIN/AEI/10.13039/501100011033 funded by the Agencia Estatal de Investigación, 
IT1615-22 funded by the Basque Government, 
and by Ikerbasque.

H. Mizutani is partially
supported by JSPS KAKENHI Grant Number JP21K03325.

L. Roncal is also supported by the projects
PID2020-113156GB-I00/MCIN/AEI/10.13039/501100011033 (acronym \lq\lq HAPDE'')
and RYC2018-025477-I 
funded by Agencia Estatal de Investigación,
and by Ikerbasque.

N. M. Schiavone is also supported by the  grant FJC2021-046835-I funded by the EU \lq\lq NextGenerationEU''/PRTR and by MCIN/AEI/10.13039/501100011033,
by the project PID2021-123034NB-I00/MCIN/ AEI/10.13039/501100011033 funded by the Agencia Estatal de Investigación,
and by the EXPRO grant No.~20-17749X of the Czech Science Foundation during his period as a researcher at Czech Technical University in Prague, where this work initiated. 
He is also member of the \lq\lq Gruppo Nazionale per L'Analisi Matematica, la Probabilit\`{a} e le loro Applicazioni'' (GNAMPA) of the \lq\lq Istituto Nazionale di Alta Matematica'' (INdAM).


\bibliographystyle{abbrv}

\begin{thebibliography}{99}
	
\bibitem{BGSZ}
F. Bagarello, J.-P. Gazeau, F. H. Szafraniec, and M. Znojil, \emph{Non-selfadjoint Operators in Quantum Physics: Mathematical aspects.} John Wiley \& Sons, 2015.

 \bibitem{BBG}
H. Bahouri, D. Barilari, and I. Gallagher, 
Strichartz estimates and Fourier restriction theorems on the Heisenberg group, \emph{J. Fourier Anal. Appl.} {\bf 27}(2) (2021).

\bibitem{BGX}
H. Bahouri, P. G\'erard, and C.-J. Xu, Espaces de Besov et estimations de Strichartz g\'en\'eralis\'ees sur le groupe de Heisenberg, \emph{J. Anal. Math.} {\bf 82} (2000), 93--118.

\bibitem{BVZ}
J.~A. Barcel{\'o}, L.~Vega, and M.~Zubeldia, 
The forward problem for the
electromagnetic {H}elmholtz equation with critical singularities, 
\emph{Adv. Math.}
\textbf{240} (2013), 636--671.
  
\bibitem{BB}
C.~M. Bender and S. Boettcher, 
Real spectra in non-Hermitian Hamiltonians having $\mathcal{PT}$ symmetry,
\emph{Phys. Rev. Lett.} {\bf 80} (1998), 5243--5246.  

 \bibitem{BPST1}
N. Burq, F. Planchon, J. Stalker, and A. S.
Tahvildar-Zadeh, Strichartz estimates for the wave and Schr\"odinger
equations with the inverse-square potential, \emph{J. Funct. Anal.} {\bf 203}
(2003), 519--549.

\bibitem{BPST2}
N. Burq, F. Planchon, J. G. Stalker, and A. S. Tahvildar-Zadeh,
Strichartz estimates for the wave and Schr\"odinger equations with potentials of critical decay,
\emph{Indiana Univ. Math. J.}, {\bf 53} (2004), 1665--1680. 

\bibitem{CCF}
B. Cassano, L. Cossetti, and L. Fanelli,
Improved Hardy--Rellich inequalities, \emph{Comm. Pure Appl. Anal.} {\bf 21} (2022), 867--889.

\bibitem{CFKP}
B. Cassano, V. Franceschi, D. Krej\v{c}i\v{r}\'ik, and D. Prandi, 
Horizontal magnetic fields and improved Hardy inequalities in the Heisenberg group, \emph{Comm. Partial Differential Equations} 
{\bf 48}(5) (2023), 711--752.

\bibitem{CFK}
	L. Cossetti, L. Fanelli and D. Krej\v{c}i\v{r}\'ik,
	Absence of eigenvalues of Dirac and Pauli Hamiltonians via the method of multipliers, 
	\emph{Comm. Math. Phys.} \textbf{379} (2020), 633--691.
	
\bibitem{CFK23}
L. Cossetti, L. Fanelli and D. Krej\v{c}i\v{r}\'ik, 
Uniform resolvent estimates and absence of eigenvalues of biharmonic operators with complex potentials. 
arXiv preprint arXiv:2309.06823 (2023).
	
\bibitem{CFS}
	L. Cossetti, L. Fanelli and N.M. Schiavone,
	Recent developments in spectral theory for non-self-adjoint Hamiltonians,
	to appear in \emph{Mathematical Physics and Its Interaction. Springer Proceedings in Mathematics \& Statistics}. 
	
\bibitem{CK}
	L. Cossetti, and D. Krej\v{c}i\v{r}\'ik,
	Absence of eigenvalues of non‐self‐adjoint Robin Laplacians on the half‐space,
	\emph{Proc. Lond. Math. Soc.} \textbf{121}.3 (2020), 584--616.
	
\bibitem{Da0}
L. D'Ambrosio, 
Critical degenerate inequalities on the Heisenberg group,
\emph{Manuscripta Mathematica}, \textbf{106}(4) (2001), 519--536.

\bibitem{Da1} 
L. D'Ambrosio, {Hardy inequalities related to Grushin type operators}, 
\emph{Proc. Amer. Math. Soc.} {\bf 132} (2004), 725--734.

\bibitem{Da2} 
L. D'Ambrosio, {Some Hardy inequalities on the Heisenberg group}, 
\emph{Differ. Equations.} {\bf 40} (2004), 552--564.

\bibitem{Da3} 
L. D'Ambrosio, {Hardy-type inequalities related to degenerate elliptic differential operators}, \emph{Ann. Sc. Norm. Super. Pisa Cl. Sci.} {\bf 5} (2005), 451--486.

\bibitem{DMW} A. Dasgupta, S. Molahajloo, and M.-W. Wong, {The spectrum of the sub-Laplacian on the Heisenberg group}, \emph{Tohoku Math. J.} {\bf 63} (2011), 269--276.

\bibitem{DFKS}
P. D’Ancona, L. Fanelli, D. Krejčiřík, and N. M. Schiavone, 
{Localization of eigenvalues for non-self-adjoint Dirac and Klein–Gordon operators}, 
\emph{Nonlinear Anal.} {\bf 214} (2022), 112565.

\bibitem{DFS}
P. D’Ancona, L. Fanelli, and N. M. Schiavone, 
Eigenvalue bounds for non-selfadjoint Dirac operators,
\emph{Math. Ann.} {\bf 383} (2022), 621--644.

\bibitem{FKV}
	L. Fanelli, D. Krej\v{c}i\v{r}\'ik, and L. Vega, 
	{Spectral stability of Schr\"odinger operators with subordinated complex potentials}, 
	\emph{J. Spectr. Theory }\textbf{8} (2018), 575--604.
	
\bibitem{FKV2}
	L. Fanelli, D. Krej\v{c}i\v{r}\'ik, and L. Vega, 
	Absence of eigenvalues of two-dimensional magnetic Schr\"odinger operators,
	\emph{J. Funct. Anal.} \textbf{275}(9) (2018), 2453--2472.

\bibitem{FMRS}
	L. Fanelli, H. Mizutani, L. Roncal, and N. M. Schiavone,
	Resolvent estimates, smoothing effects and spectral stability for the Heisenberg sublaplacian,
	in preparation.


\bibitem{F} G. B. Folland,
A fundamental solution for a subelliptic operator,
\emph{Bull. Amer. Math. Soc.} {\bf 79} (1973),
373--376.

 \bibitem{F2} G. B. Folland, \emph{Harmonic Analysis in Phase Space}, Ann. Math. Stud. {\bf 122}. Princeton University Press, Princeton, N.J., 1989.
 
 \bibitem{FS} G. B. Folland and E. M. Stein,
{Estimates for the $\bar{\partial}_b$ complex and analysis on the Heisenberg group,}
\emph{Comm. Pure Appl. Math.} {\bf 27} (1974),
429--522.

 
 \bibitem{Frank}
	R. L. Frank, 
	Eigenvalue bounds for Schr\"odinger operators with complex potentials, 
	\emph{Bull. Lond. Math. Soc.} \textbf{43} (2011), 745--750.

\bibitem{FrankIII}
	R. L. Frank, 
	Eigenvalue bounds for Schr\"odinger operators with complex potentials. III,
	\emph{Trans. Amer. Math. Soc.} \textbf{370}(1) (2018), 219--240.
	
\bibitem{FL}
R. L. Frank and E. Lieb, 
Sharp constants in several inequalities on the Heisenberg group, \emph{Annals of Math.} {\bf 176} (2012), 349--381.	
	
\bibitem{FS17}
	R. L. Frank and B. Simon,
	Eigenvalue bounds for Schr\"odinger operators with complex potentials. II,
	\emph{J. Spectr. Theory} \textbf{7} (2017), 633--658. 
 
\bibitem{GL} N. Garofalo and E. Lanconelli, {Frequency functions on the Heisenberg group, the uncertainty principle and unique continuation}, \emph{Ann. Inst. Fourier (Grenoble)} {\bf 40} (1990), 313--356.
 
\bibitem{Greiner}
	P. C. Greiner,
	Spherical harmonics on the Heisenberg group,
	\emph{Canad. Math. Bull.} \textbf{23}(4) (1980), 383--396. 
	
\bibitem{GK}
	P. C. Greiner and T. H. Koornwinder,
	Variations on the Heisenberg spherical harmonics,
	\emph{Stichting Mathematisch Centrum. Zuivere Wiskunde} ZW 186/83 (1983).

\bibitem{Gut}
	S. Gutiérrez, 
	Non trivial $L^q$ solutions to the Ginzburg-Landau equation,
	\emph{Math. Ann.} \textbf{328} (2004), 1--25.
	
\bibitem{HK}
	M. Hansmann, and D. Krej\v{c}i\v{r}\'ik, 
	The abstract Birman--Schwinger principle and spectral stability. 
	\emph{J. Anal. Math.} \textbf{148}(1) (2022), 361--398.
 

\bibitem{IS}
T.~Ikebe and Y.~Saito, 
Limiting absorption method and absolute continuity
  for the {S}chr{\"o}dinger operator, \emph{J. Math. Kyoto Univ.} \textbf{12} (1972),
  513--542.
  
  \bibitem{JL}
D.~Jerison and J. M.~Lee, 
Extremals for the Sobolev inequality on the Heisenberg group and the CR Yamabe problem, \emph{J. Amer. Math. Soc.} \textbf{1} (1988),
  1--13.

\bibitem{K}
T.~Kato,
 \emph{Perturbation theory for linear operators}, Springer-Verlag,
  Berlin, 1966.
  
\bibitem{KY}
T.~Kato and K. Yajima,
Some examples of smooth operators and the associated smoothing effect, 
\emph{Rev. Math. Phys.} {\bf 1}(4) (1989), 481--496.  

\bibitem{KRS}
 C. E. Kenig, A. Ruiz, and C. D. Sogge, Uniform Sobolev inequalities and unique continuation for second order constant coefficient differential operators, \emph{Duke Math. J.} {\bf 55} (1987), 329--347.
 
\bibitem{KL}
	Y. Kwon and S. Lee, 
	Sharp resolvent estimates outside of the uniform boundedness range,
	\emph{Comm. Math. Phys.} \textbf{374}(3) (2020), 1417--1467.
 

\bibitem{MS}
H. Mizutani, and N. M. Schiavone, Spectral enclosures for Dirac operators perturbed by rigid potentials, \emph{Rev. Math. Phys.} \textbf{34}(8) (2022), 2250023.

\bibitem{ReedSimonIV}
M. Reed, and B. Simon,
\emph{Methods of Modern Mathematical Physics IV: Analysis of Operators.} 
Academic Press Inc., 1978.

\bibitem{R}
F.~Rellich,
Halbbeschr\"{a}nkte Differentialoperatoren h\"{o}herer Ordnung. 
In \emph{Proceedings of the International Congress of Mathematicians}, 1954, Amsterdam, vol. III, 243--250. Erven P. Noordhoof N.V., Groningen, 1956. 

\bibitem{RXZ}
	T. Ren, Y. Xi, and C. Zhang,
	An endpoint version of uniform Sobolev inequalities,
	\emph{Forum Math.} \textbf{30}(5) (2018), 1279--1289.

\bibitem{RS}
M. Ruzhansky, and D. Suragan,
\emph{Hardy Inequalities on Homogeneous Groups: 100 years of Hardy inequalities.} Springer Nature, 2019.

  
\bibitem{SGY}
F.G. Scholtz, H.B. Geyer, and F.J.W. Hahne,
Quasi-Hermitian operators in quantum mechanics and the variational principle, \emph{Ann. Phys.} {\bf 213} (1992), 74--101.  
  
 \bibitem{S}
B.~Simon, \emph{Quantum Mechanics for {H}amiltonians Defined as Quadratic
  Forms}, Princeton Univ. Press., New Jersey, 1971.


\end{thebibliography}

\end{document}